\documentclass[11pt,reqno]{amsart}

\newcommand{\mysection}[1]{\section{#1}
      \setcounter{equation}{0}}

\newtheorem{theorem}{Theorem}[section]
\newtheorem{lemma}[theorem]{Lemma}

\newtheorem{corollary}[theorem]{Corollary}

\theoremstyle{definition}
\newtheorem{assumption}{Assumption}[section]
\newtheorem{definition}{Definition}[section]
\theoremstyle{remark}
\newtheorem{remark}{Remark}[section]

\newcommand\cbrk{\text{$]$\kern-.15em$]$}} 
\newcommand\opar{\text{\raise.2ex\hbox{${\scriptstyle | }$}\kern-.34em$($} }

 \makeatletter
 \def\dashint{%
 \operatorname%
 {\,\,\text{\bf--}\kern-.98em\DOTSI\intop\ilimits@\!\!}}
 \makeatother

\newcommand\gb{\mathfrak{b}}

\newcommand\bR{\mathbb{R}}

\newcommand\bL{\mathbb{L}}

\newcommand\bW{\mathbb{W}}

\newcommand\cD{\mathcal{D}}

\newcommand\cL{\mathcal{L}}

\newcommand\cW{\mathcal{W}}

\newcommand\frD{\mathfrak{D}}

\begin{document}

\title[Divergence equations
with growing   coefficients]{
On divergence form second-order PDEs
with growing  coefficients in $W^{1}_{p}$
spaces without weights}

\author[N.V.  Krylov]{N.V. Krylov}%
\thanks{The work   was partially supported
  by NSF grant DMS-0653121}
\address{127 Vincent Hall, University of Minnesota,
Minneapolis,
       MN, 55455, USA}
\email{krylov@math.umn.edu}


\begin{abstract}
We consider second-order divergence form uniformly parabolic 
and elliptic PDEs with
bounded  and $VMO_{x}$  leading coefficients and
possibly linearly growing lower-order coefficients.
 We look for solutions
  which are summable to the $p$th power with respect to the usual
Lebesgue measure along with their first   derivatives
with respect to the spatial variables.
\end{abstract}

\maketitle

\mysection{Introduction}

We consider divergence form uniformly parabolic and elliptic
second-order  PDEs with
bounded  and $VMO_{x}$ leading coefficients and
possibly linearly growing lower-order coefficients.
 We look for solutions
  which are summable to the $p$th power with respect to the usual
Lebesgue measure along with their first   derivatives
with respect to the spatial variables.  
In some sense we extend the results of \cite{Kr10_1},
where $p=2$,
to general $p\in(1,\infty)$. However in \cite{Kr10_1}
there is no regularity assumption on the leading coefficients
and there are also stochastic terms in the equations.

As in \cite{CV} one of the main motivations for studying
 PDEs with growing first-order coefficients
 is filtering theory 
for partially observable diffusion processes.
 
 It is generally believed that introducing weights is the 
most natural setting for equations with growing 
coefficients. When the 
coefficients grow it is quite natural to consider the equations
in function 
 spaces with weights that would restrict the set of solutions
in such a way that all terms in the equation will be
from the same space as the free terms. 
 The present paper seems to be the first
one treating the unique solvability of
these equations with growing
lower-order coefficients in the usual Sobolev spaces $W^{1}_{p}$
 without weights and
without imposing any {\em special\/}
 conditions on the relations between
the coefficients or on their {\em derivatives\/}.  

The theory of  PDEs and stochastic PDEs in Sobolev  spaces {\em
with\/} weights
attracted some attention in the past. We do not use weights and only
mention  a few papers about  stochastic 
PDEs in $\cL_{p}$-spaces
 with weights
in which one can find further references: \cite{AM} (mild solutions,
general $p$),
\cite{CV}, \cite{Gy93}, \cite{Gy97}, \cite{GK} ($p=2$ in the four
last articles).

Many more papers are devoted to the theory of {\em deterministic\/}
 PDEs
with growing coefficients in Sobolev spaces with weights. 
We  cite only a few of them sending the reader to the references
therein again because neither  do we deal with weights nor  use
the results of these papers. It is also worth saying that our results
do not generalize the results of these papers.

In most of them the coefficients
are time independent, see \cite{CV1}, \cite{ChG}, \cite{FL},  
 \cite{MP1}, part of the result of which
are extended in \cite{GL} to
time-dependent Ornstein-Uhlenbeck operators.

It is worth noting that many issues for  deterministic 
divergence-type equations with time independent 
growing coefficients in $\cL_{p}$ spaces with arbitrary $p
\in(1,\infty)$ {\em without\/} weights
were also treated previously in the literature. 
This was done mostly  by using the semigroup approach
which excludes time dependent coefficients
and makes it almost impossible to use the results in
the more or less general filtering theory.
We briefly mention only a few recent papers
sending the reader to them for additional information.

In \cite{LV} a strongly continuous
in $\cL_{p}$ semigroup is constructed corresponding
to elliptic  operators with measurable
leading coefficients and Lipschitz
continuous drift coefficients.  
In \cite{MP} it is assumed that
 if, for $| x|\to\infty$, the drift coefficients   grow,
then  the zeroth-order coefficient  should grow, basically,
as the square of the drift. There  is also a condition on the divergence
of the drift coefficient.
In \cite{PR} there is no zeroth-order term
and the semigroup is constructed under some assumptions
one of which translates into  the monotonicity of 
$\pm b(x)-Kx$, for a constant $K$, if the leading term is the Laplacian.
In \cite{CF}  the drift coefficient
is assumed to be globally Lipschitz
continuous if the zeroth-order coefficient is constant.

Some conclusions in the above cited papers are quite similar to ours
but the corresponding assumptions are not as general
in what concerns the regularity of the coefficients.
However, these papers contain a lot of additional
important information not touched upon in the present paper
(in particular, it is shown in \cite{LV} that the corresponding semigroup
is not analytic and in \cite{Me} that the spectrum of
an elliptic operator in $\cL_{p}$ {\em depends\/} on $p$).

The technique, we apply, originated from \cite{KP}
and \cite{Kr09_1} and uses special
 cut-off functions whose support evolves in time
in a manner adapted to the drift. As there,
we do not make any regularity assumptions on the coefficients
in the time variable but unlike \cite{Kr10_1}, where $p=2$, we use 
 the results of  \cite{Kr07} where  
some regularity on the coefficients in $x$ variable
is needed, like, say, the condition that the 
second order coefficients
be in VMO uniformly with respect to the time variable.
   
It is worth noting that considering
divergence form equations in $\cL_{p}$-spaces
is quite useful in the treatment of
filtering problems (see, for instance, \cite{Kr_10}) especially
when the power of summability is taken large  and we
intend to treat this issue in a subsequent paper.
 
The article is organized as follows. In Section
\ref{section 6.9.1} we describe the problem,
Section \ref{section 2.15.1} contains the statements
of two main results, Theorem \ref{theorem 3.11.1}
on an apriori estimate  providing,
in particular, uniqueness of solutions
 and Theorem \ref{theorem 3.16.1}
about the existence of solutions. The results about Cauchy's problem
and elliptic equations are also given there.
Theorem \ref{theorem 3.11.1} is proved in 
Section~\ref{section 6.9.3}
after we prepare the necessary tools in Section
\ref{section 6.9.2}.
Theorem \ref{theorem 3.16.1} is proved in the last
Section \ref{section 6.9.5}.

As usual when we speak of
 ``a constant" we always mean ``a finite constant".

 The author discussed the article with Hongjie Dong
 whose comments are greatly appreciated.

\mysection{Setting of the problem}
                                       \label{section 6.9.1}

 We consider the   second-order
operator $L_{t}$
$$
 L_{t} u_{t}(x) =D_{i}\big( a^{ij}_{t}( x)D_{j}u_{t}(x)
+\gb^{i}_{t}(x)u_{t}(x)\big) + b^{i} _{t}(x)
D_{i} u_{t}(x) -c _{t}(x) u_{t}(x),
$$
   acting on functions $u_{t}(x)$ defined
on $([S,T]\cap\bR) \times \bR^d$ 
(the summation convention is enforced throughout the article),
where $S$ and $T$ are   such that $-\infty\leq S
<T\leq\infty$. Naturally,
$$
D_{i}=\frac{\partial}{\partial x^{i}}
$$

Our main concern is proving the unique solvability
of the equation
\begin{equation}
                                                \label{2.6.4}
\partial_{t}u_{t}=
L_{t}u_{t}-\lambda u_{t}+D_{i}f^{i}_{t}+f^{0}_{t}
\quad t\in[S,T]\cap\bR,
\end{equation}
with an appropriate initial condition at $t=S$ if $S>
-\infty$, where  
$\lambda>0$ is a constant and $\partial_{t}=\partial/\partial t$. 
The precise 
assumptions on the coefficients, free terms, and initial data
will be given later. First we introduce appropriate function
spaces.

Denote $C^{\infty}_{0}=C^{\infty}_{0}(\bR^{d})$,
$\cL_{p}=\cL_{p}(\bR^{d})$, and let $W^{1}_{p} =W^{1}_{p}(\bR^{d})$
be the Sobolev space of functions $u$ of class
 $\cL_{p}$, such that
$Du\in \cL_{p}$, where $Du$ is the gradient of $u$
and $1<p<\infty$. 
For $-\infty\leq S<T\leq\infty$ define
$$
\bL_{p}(S,T)=\cL_{p}((S,T),\cL_{p} ),
\quad
\bW^{1}_{p}(S,T)=\cL_{p}((S,T),
W^{1}_{p} ),
$$
$$
\bL_{p}(T)=\bL_{p}(-\infty,T),\quad
\bW^{1}_{p}(T)=\bW^{1}_{p}(-\infty,T),
$$
$$
\bL_{p}=\bL_{p}(\infty),\quad
\bW^{1}_{p} =\bW^{1}_{p}( \infty ).
$$
Remember that the elements of
$\bL_{p}(S,T)$ need only belong to $\cL_{p}$
on a  Borel subset of 
$(S,T)$ of full measure. We will always assume
that these elements  are defined everywhere
on $(S,T)$ at least as generalized functions
on $\bR^{d}$.
Similar situation occurs in the case of $\bW^{1}_{p}(S,T)$.
 
The following definition is most appropriate
for investigating our equations
if the coefficients of $L$ are bounded.
 
 \begin{definition}
                                         \label{definition 3.16.1}
     
We introduce the space $\cW^{1}_{p}(S,T)$,
which is the space of functions $u_{t}
$ on $[S,T]\cap\bR $ with values
in the space of generalized functions on $\bR^{d}$
and having the following properties:

(i)  We have $u
\in \bW^{1}_{p}(S,T)$;

(ii) There exist   $f^{i}\in \bL_{p}(S,T)$,
$i=0,...,d$,  
such that
 for any $\phi\in C^{\infty}_{0}$ and finite $s,t\in[S,T] $
we have
\begin{equation}
                                                 \label{1.2.1}
(u_{t},\phi)=(u_{s},\phi)+
 \int_{s}^{t}\big(
 (f^{0}_{r},\phi)-(f^{i}_{r},D_{i}\phi)\big)\,dr.
\end{equation}
In particular, for any $\phi\in C^{\infty}_{0}$, the function
$(u_{t},\phi)$ is   continuous on $[S,T]\cap\bR$.
In case that property (ii) holds, we write
$$
\partial_{t}u_{t}=D_{i}f^{i}_{t}+f^{0}_{t} ,
\quad t\in[S,T]\cap\bR.
$$
 
\end{definition}

Definition \ref{definition 3.16.1} allows us to introduce
the spaces of initial data
\begin{definition}
                                        \label{definition 8.10.1}
Let $g$ be a generalized function. We write
$g\in W^{1-2/p}_{p}$ if there exists a function $v_{t}\in \cW^{1}_{p}
(0,1)$ such that $\partial_{t}v_{t}=\Delta v_{t}$, $t\in[0,1]$,
and $v_{0}=g$. In such a case we set
$$
\|g\|_{W^{1-2/p}_{p}}=\|v\|_{\bW^{1}_{p}(0,1)}.
$$
\end{definition}

Following Definition \ref{definition 3.16.1}
we understand
equation \eqref{2.6.4} as the requirement that 
for any  $\phi\in C^{\infty}_{0}$ and
finite $s,t\in[S,T]$ we have
$$
(u_{t},\phi)=(u_{s},\phi)
+
\int_{s}^{t}\big[(b^{i}_{r}D_{i}u_{r}
-(c_{r}+\lambda)u_{r}+f^{0}_{r},\phi)
$$
\begin{equation}
                                                \label{3.16.7}
-(a^{ij}_{r}D_{j}u_{r}+\gb^{i}_{r}u_{r}+
f^{i}_{r},D_{i}\phi)
 \big]\,dr.
\end{equation}

Observe that at this moment it is not clear that the right-hand
side makes sense.
Also notice that, if the coefficients of $L$  
are bounded, then any $u\in\cW^{1}_{p}(S,T)$ is a solution
of \eqref{2.6.4} with appropriate free terms since if
\eqref{1.2.1} holds, then \eqref{2.6.4} holds as well with
$$
f^{i}_{t}-a^{ij}_{t}D_{j}u_{t}-\gb^{i}_{t}u_{t},
\quad i=1,...,d,\quad
f^{0}_{t}+(c_{t}+\lambda)u_{t}-b^{i}_{t}D_{i}u_{t},
$$
in place of $f^{i}_{t}$, $i=1,...,d$, and $f^{0}_{t}$,
respectively.

 We  give the definition of
solution of  \eqref{2.6.4} adopted throughout the article
and which in case the coefficients of $L$
are bounded coincides with the one obtained by applying 
 Definition \ref{definition 3.16.1}.
 
\begin{definition}
                                      \label{definition 3.20.01}
Let $f^{j}\in\bL_{p}(S,T)$, $j=0,...,d$ and assume that $S>-\infty$.
By a solution of
\eqref{2.6.4}   with initial condition
$u_{S}\in W^{1-2/p}_{p}$
we mean a function $u\in\bW^{1}_{p}(S,T)$ (not
$\cW^{1}_{p}(S,T)$) such that

(i)  For any $\phi\in C^{\infty}_{0} $ 
 the integral
with respect to $dr$ in \eqref{3.16.7} is 
well defined and is finite for all finite
 $s,t\in[S,T]$;

(ii) For any $\phi\in C^{\infty}_{0} $
equation  \eqref{3.16.7}  
 holds for all finite
 $s,t\in[S,T]$.

In case $S=-\infty$ we drop mentioning initial 
condition in the above lines.

\end{definition}

\mysection{Main results}
                                          \label{section 2.15.1}

For $\rho>0$ denote $B_{\rho}(x)=\{y\in\bR^{d}:|x-y|<\rho\}$,
$B_{\rho}=B_{\rho}(0)$. 
\begin{assumption} 
                                       \label{assumption 2.7.2}
(i) The functions $a^{ij}_{t}(x)$, $\gb^{i}_{t}(x)$,
$b^{i}_{t}(x)$,  and $c_{t}(x)$  are real valued and Borel measurable 
 and  $c \geq 0$.

(ii) There exists  a constant $\delta>0$ such that
for all values of arguments and $\xi\in\bR^{d}$
$$
a^{ij}  \xi^{i}
\xi^{j}\geq\delta|\xi|^{2},\quad
|a^{ij}|\leq \delta^{-1}  .
$$
 Also, the constant $\lambda>0$.

(iii) For any  $x\in
\bR^{d}$  the function
$$
\int_{B_{1}}(|\gb_{t}(x+y)|+|b_{t}(x+y)|+c_{t}(x+y) )\,dy
$$ 
is locally integrable to the $p'$th power on $\bR$,
where $p'=p/(p-1)$.
\end{assumption}

Notice that the matrix $a=(a^{ij})$ need not be symmetric.
Also notice that in Assumption \ref{assumption 2.7.2} (iii)
the ball $B_{1}$ can be replaced with any other ball
without changing the set of admissible coefficients
$\gb,b,c$.

We take and fix constants $K\geq0,\rho_{0},\rho_{1}
\in(0,1]$, and choose
a number $q=q(d,p) $ so  that
\begin{equation}
                                                  \label{8.11.3}
q>\min(d,p ),\quad q>\min(d,p' ),\quad q\geq\max(d, p,p').
\end{equation}

The following assumptions contain   a  parameter  $ 
\gamma \in(0,1]$,  
whose value will be specified later.

 \begin{assumption} 
                                      \label{assumption 3.11.1} 
For   $ \gb:=(\gb^{1} ,...,\gb^{d} )$
and $ b:=(b^{1} ,...,b^{d} )$ and   $(t,x)\in\bR^{d+1}$
 we have  
$$
 \int_{B_{\rho_{1}}(x)}\int_{B_{\rho_{1}}(x)}|\gb_{t}( y)
-\gb_{t}(   z)|^{q }\,dydz  +
 \int_{B_{\rho_{1}}(x)}\int_{B_{\rho_{1}}(x)}|b_{t}( y)
-b_{t}(  z)|^{q}\,dydz   
$$
$$
+ \int_{B_{\rho_{1}}(x)}\int_{B_{\rho_{1}}(x)}|c_{t}( y)
-c_{t}( z)|^{q}\,dydz  \leq KI_{q>d}+\rho_{1}^{d}\gamma.
$$

\end{assumption}

Obviously,  Assumption \ref{assumption 3.11.1}  
is satisfied   if  
 $b$, $\gb$, and $c$  are independent of $x$.
They also are satisfied  with any $q>d$, $\gamma\in(0,1]$,
and $\rho_{1}=1$ on the 
account of choosing
$K$ appropriately if, say,
$$
|\gb_{t}(  x )
-\gb_{t}(y)|+|b_{t}(  x )
-b_{t}( y)|+|c_{t}(  x )
-c_{t}(  y)|\leq N
$$ 
whenever $|x-y|\leq1$, where $N$ is a constant.
 We see   that Assumption \ref{assumption 3.11.1}
 allows $b$, $\gb$, and $c$ growing linearly in $x$.

 \begin{assumption} 
                                       \label{assumption 2.4.01}
For any  
$\rho\in(0,\rho_{0}]$, $s\in\bR $, 
 and $i,j=1,...,d$  we have 
\begin{equation}
                                             \label{2.9.5}
\rho^{-2d-2}\int_{s}^{s+\rho^{2}}\bigg(
 \sup_{x\in \bR^{d}}\int_{B_{\rho}(x)}\int_{B_{\rho}(x)}
|a^{ij}_{t}( y)-a^{ij}_{t}( z)|\,dydz\bigg)
 \,dt\leq\gamma .
\end{equation}

\end{assumption}

Obviously, the left-hand side of \eqref{2.9.5}
is less than
$$
N(d)\sup_{t \in\bR}\sup_{|x-y|\leq2\rho }|a^{ij}_{t}( x )
-a^{ij}_{t}(y)|,
$$
which implies that Assumption \ref{assumption 2.4.01}
 is satisfied with any $\gamma \in(0,1]$ if, for instance,
$a$ is   uniformly continuous in $x$ uniformly with respect to $t$.
Recall that if $a$ is independent of $t$ and
 for any $\gamma >0$ there is a $\rho_{0}>0$
such that Assumption~\ref{assumption 2.4.01}
 is satisfied, then one says that
$a$ is in VMO.

\begin{theorem}
                                       \label{theorem 3.11.1}
There exist
$$
\gamma =\gamma (d,\delta,p  )\in(0,1],
$$
$$
 N=N(d,\delta,p   ),
\quad \lambda_{0}=\lambda_{0}(d,\delta,  p,
\rho_{0},\rho_{1},K)\geq1
$$  
such that, if the above assumptions are satisfied
and $\lambda\geq
\lambda_{0}$ and
  $u  $ is a solution of \eqref{2.6.4}
with  zero initial data $($if $S>-\infty$$)$ and some
$f^{j} \in\bL_{p}(S,T)$, then
\begin{equation}
                                       \label{3.11.2}
\lambda\|u\|^{2}_{\bL _{p}(S,T)}+\|Du\|^{2}_{\bL _{p}(S,T)}
\leq N\big(\sum_{i=1}^{d}\|f^{i}\|^{2}_{\bL _{p}(S,T)}
+\lambda^{-1}\|f^{0}\|^{2}_{\bL _{p}(S,T)} \big).
\end{equation}
\end{theorem}

Notice that the main case of Theorem \ref{theorem 3.11.1}
is when $S=-\infty$ because if $S>-\infty$ and $u_{S}=0$,
then the function $u_{t}I_{t\geq S}$ will be a solution
of our equation on $(-\infty,T]\cap\bR$ with
$f^{j}_{t}=0$ for $t<S$.

This theorem provides an apriori estimate implying
uniqueness of solutions.
Observe that the assumption that such a solution
exists is quite nontrivial because   if $\gb_{t}(x)\equiv x$,
it is not true that  $\gb u\in\bL_{p}(S,T)$
for arbitrary $u\in\bW^{ 1}_{p}(S,T)$.

It is also worth noting that, as can be easily 
seen from
the proof of Theorem \ref{theorem 3.11.1}, one can choose
a function $\gamma(d,\delta,p)$ so that it is
continuous  in $(\delta,p)$. The same holds for $N$ and $\lambda_{0}$
from Theorem \ref{theorem 3.11.1}.

We have a similar result for nonzero initial data.

\begin{theorem}
                                       \label{theorem 8.8.1}
Let $S>-\infty$. In Theorem \ref{theorem 3.11.1}
replace the assumption that $u_{S}=0$ with the assumption
that $u_{S}\in W^{1-2/p}_{p}$. Then its statement
remains true if in the right-hand side of \eqref{3.11.2}
we add the term
$$
N \|u_{S}\|^{2}_{W^{1-2/p}_{p}}.
$$
\end{theorem}

Proof. Take $v_{t}$ from Definition \ref{definition 8.10.1}
corresponding to $g=u_{S}$
and
set
$$
\tilde{u}_{t}=\begin{cases}
u_{t}& t\geq S,\\
(t-S+1)v_{S-t}& S\geq t\geq  S-1,\\
0& S-1\geq t 
\end{cases}
$$
and for $i=1,...,d$ set 
$$
\tilde{f}^{i}_{t}=\begin{cases}
f^{i}_{t}& t\geq S,\\
 -2(t-S+1)D_{i}v_{S-t}& S> t\geq  S-1,\\
0& S-1\geq t ,
\end{cases}
$$
$$
\tilde{f}^{0}_{t}=\begin{cases}
f^{0}_{t}& t\geq S,\\ \,
 [1+\lambda(t-S+1)]v_{S-t}& S> t\geq  S-1,\\
0& S-1\geq t .
\end{cases}
$$
We also modify the coefficients of $L$ by 
multiplying each one of them but $a^{ij}_{t}$
by $I_{t\geq S}$ and setting
$$
\tilde{a}^{ij}_{t}=\begin{cases}
a^{ij}_{t}& t\geq S,\\
\delta^{ij}& S>t.
\end{cases}
$$
Here we profit from the fact that no regularity assumption
on the dependence of the coefficients on $t$ is imposed.
By denoting by $\tilde{L}$ the operator with the
modified coefficients we easily see that
$\tilde u_{t}$ is a solution (always in the sense
of Definition \ref{definition 3.20.01})
of
$$
\partial_{t}\tilde u_{t}= \tilde
L_{t}\tilde u_{t}- \lambda \tilde u_{t}+D_{i}
\tilde f^{i}_{t}+\tilde f^{0}_{t},
\quad t\leq T.
$$

By Theorem \ref{theorem 3.11.1}
$$
\lambda\| u\|^{2}_{\bL _{p}(S,T)}+
\|D u\|^{2}_{\bL _{p}(S,T)}
\leq N\big(\sum_{i=1}^{d}\|\tilde f^{i}\|^{2}_{\bL _{p}( T)}
+\lambda^{-1}\|\tilde f^{0}\|^{2}_{\bL _{p}( T)} \big),
$$
where
$$
\|\tilde f^{i}\|^{p}_{\bL _{p}( T)}=
\| f^{i}\|^{p}_{\bL _{p}(S, T)}
+\|\tilde f^{i}\|^{p}_{\bL _{p}(S-1, S)}\leq
\|  f^{i}\|^{p}_{\bL _{p}(S, T)}+2^{p}
\|D_{i}v\|^{p}_{\bL_{p}(0,1)} 
$$
$$
\leq
\| f^{i}\|^{p}_{\bL _{p}(S, T)}+2^{p}
 \|u_{S}\|_{W^{1-2/p}_{p}}^{p},
$$
$$
\|\tilde f^{0}\|^{p}_{\bL _{p}( T)} 
\leq \|f^{0}\|^{p}_{\bL _{p}(S,T)}+N(1+\lambda^{p})
\|v\|_{\bL _{p}( 0,1)}^{p}
$$
$$
\leq \|f^{0}\|^{p}_{\bL _{p}(S,T)}+N 
(1+\lambda^{p})\|u_{S}\|_{W^{1-2/p}_{p}}^{p}.
$$
Since $\lambda\geq\lambda_{0}\geq1$, we have $1+\lambda^{p}
\leq 2\lambda^{p}$ and we
 get our assertion
thus proving the theorem.

  Here is an existence theorem.

\begin{theorem}
                                        \label{theorem 3.16.1}
Let the above assumptions be satisfied with
$\gamma $ taken from Theorem \ref{theorem 3.11.1}.
Take   $\lambda\geq\lambda_{0}$, where $\lambda_{0}$
is defined in Theorem \ref{theorem 3.11.1}.  Then
for any $f^{j}\in\bL_{p}(T)$, $j=0,...,d$,
 there exists
a unique solution of \eqref{2.6.4} with $S=-\infty$.
\end{theorem}

 It turns out that the solution, if it exists, is independent
of the space in which we are looking for solutions.

\begin{theorem}
                                               \label{theorem 9.26.1}
Let $ 1<p_{1}\leq p_{2}<\infty$ and let
$$
\gamma=\inf_{p\in[p_{1},p_{2}] }\gamma(d,\delta,p),
$$
where $\gamma(d,\delta,p)$ is taken from Theorem \ref{theorem 3.11.1}.
Suppose that Assumptions \ref{assumption 2.7.2}
 through \ref{assumption 2.4.01}
are satisfied with so defined $\gamma$ and with $p=p_{1}$ and $p=p_{2}$.

(i) Let $-\infty<S<T<\infty$,
$f^{j}\in\bL_{p_{1}} (S,T)\cap \bL_{p_{2}}(S,T) $, $j=0,...,d$,
$u_{S}\in W^{1-2/p_{1}}_{p_{1}}\cap W^{1-2/p_{2}}_{p_{2}}$,
and let $u\in\bW^{1}_{p_{1}}(S,T) \cup\bW^{1}_{p_{2}}(S,T) $ be a solution of
\eqref{2.6.4}.
  Then 
$u\in\bW^{1}_{p_{1}}(S,T)\cap\bW^{1}_{p_{2}}(S,T)$.

(ii) Let $S=-\infty, T=\infty$,
$f^{j}\in\bL_{p_{1}} \cap \bL_{p_{2}} $, $j=0,...,d$,
and let $u\in\bW^{1}_{p_{1}} \cup\bW^{1}_{p_{2}} $ be a solution of
\eqref{2.6.4}  with
\begin{equation}
                                                      \label{9.27.1}
\lambda\geq\sup_{p\in[p_{1},p_{2}] }\lambda_{0}(d,\delta,  p,
\rho_{0},\rho_{1},K),
\end{equation}
where $\lambda_{0}(d,\delta,  p,
\rho_{0},\rho_{1},K)$ is taken from Theorem 
\ref{theorem 3.11.1}. Then 
$u\in\bW^{1}_{p_{1}} \cap\bW^{1}_{p_{2}} $.
\end{theorem}

This theorem is proved in Section \ref{section 6.9.5}.
The following theorem is about Cauchy's problem with
nonzero initial data.
\begin{theorem}
                                        \label{theorem 8.8.2}
Let $S>-\infty$ and take a function $u_{S}
\in W^{1-2/p}_{p}$.
Let the above assumptions be satisfied with
$\gamma $ taken from Theorem \ref{theorem 3.11.1}.
Take   $\lambda\geq\lambda_{0}$, where $\lambda_{0}$
is defined in Theorem \ref{theorem 3.11.1}.  Then
for any $f^{j}\in\bL_{p}(S,T)$, $j=0,...,d$,
 there exists
a unique solution of \eqref{2.6.4} with 
initial value $u_{S}$.
\end{theorem}
Proof. As in the proof of Theorem \ref{theorem 8.8.1} we extend
our coefficients and $f^{j}_{t}$ for $t<S$ and then find
a unique solution $\tilde{u}_{t}$ of
$$
\partial_{t}\tilde u_{t}=
\tilde{L}_{t}\tilde{u}_{t}- \lambda \tilde{u}_{t}
+D_{i}\tilde{f}^{i}_{t}+
\tilde{f}^{0}_{t}
\quad t\in(-\infty, T]\cap\bR,
$$
By construction $(t-S+1)v_{S-t} $ satisfies this equation
for $t\leq S$, so that by uniqueness (Theorem \ref{theorem 3.11.1}
with $S$ in place of $T$) it coincides with
$\tilde{u}_{t}$ for $t\leq S$. In particular, $\tilde{u}_{S}=
v_{0}=u_{S}$. Furthermore $\tilde{u}$ satisfies \eqref{2.6.4}
since the coefficients of $\tilde{L}_{t}$
coincide with the corresponding coefficients of $L_{t}$
for finite $t\in[S,T]$. The theorem is proved.

\begin{remark}
If both $S$ and $T$ are finite, then in the above theorem
one can take $\lambda=0$. To show this take a large $\lambda>0$
and replace
 the unknown function $u_{t}$ with $v_{t}
e^{ \lambda t}$. This leads to an equation for $v_{t}$
with the additional term $-\lambda v_{t} $ and the free terms
multiplied by $e^{-\lambda t}$. The existence
of solution $v$ will be then equivalent
to the existence of $u$ if $S$ and $T$ are finite.
                                                              
\end{remark}
\begin{remark}
From the above proof and from Theorem \ref{theorem 9.26.1}
it follows that the solution, if it exists, is independent
of $p$ in the same sense as in Theorem~\ref{theorem 9.26.1}.

\end{remark}
Here is a result for elliptic equations.

\begin{theorem}
                                               \label{theorem 8.8.3}
Let the coefficients of $L_{t}$ be independent
of $t$, so that we can set $L=L_{t}$ and drop the subscript
$t$ elsewhere, let
Assumptions \ref{assumption 2.7.2} (i), (ii)  
 be satisfied,
and let $\gb$, $b$, and $c$ be locally integrable.
Then
there exist
$$
\gamma =\gamma (d,\delta,p  )\in(0,1],
$$
$$
 N=N(d,\delta,p   ),
\quad \lambda_{0}=\lambda_{0}(d,\delta, p,
\rho_{0},\rho_{1},K)\geq1
$$  
such that, if Assumptions \ref{assumption 3.11.1} and    
\ref{assumption 2.4.01} are satisfied
and $\lambda\geq
\lambda_{0}$ and
  $u  $ is a $W^{1}_{p}$-solution of 
\begin{equation}
                                                  \label{8.8.6}
Lu- \lambda  u+D_{i}f^{i}+f^{0}=0
\end{equation}
in $\bR^{d}$ with some $f^{j}\in\cL_{p}$, $j=0,...,d$,
 then
\begin{equation}
                                               \label{8.12.2}
\lambda\|u\|^{2}_{\cL _{p}}+\|Du\|^{2}_{ \cL _{p}}
\leq N\big(\sum_{i=1}^{d}\|f^{i}\|^{2}_{ \cL _{p}}
+\lambda^{-1}\|f^{0}\|^{2}_{ \cL _{p}} \big).
\end{equation}
 
Furthermore, for any $f^{j}\in\cL_{p}$, $j=0,...,d$, 
and $\lambda\geq
\lambda_{0}$ there exists
a unique solution $u\in W^{1}_{p}$ of \eqref{8.8.6}.

\end{theorem}

This result is obtained from the previous ones
in a standard way (see, for instance,
the proof of Theorem 2.1 of \cite{Kr09_1}).
One of remarkable features of \eqref{8.12.2} is that $N$
is independent of $\gb$, $b$, and $c$. It is remarkable
even if they are constant, when there is no 
assumptions on them apart from $c\geq0$.
Another point worth noting is that if $\gb=b\equiv0$,
then for the solution $u$ we have $cu\in W^{-1}_{p }$.
However, generally it is not true that 
$cu\in W^{-1}_{p }$ for {\em any} $u\in W^{1}_{p}$. For 
instance $u(x):=(1+|x|)^{-1}\in W^{1}_{p}$ if $p>d$, but if $c(x)=|x|$,
then $(1-\Delta)^{-1/2}(cu)(x)\to 1$ as $|x|\to\infty$ and
$(1-\Delta)^{-1/2}(cu)$ is not integrable to any power $r>1$.
Therefore generally,   $(L-\lambda)W^{1}_{p}\supset W^{-1}_{p}$
with proper inclusion, that does not happen if the coefficients
of $L$ are bounded.

\begin{remark}
It follows, from the arguments leading to the proof 
of Theorem \ref{theorem 8.8.3}(see \cite{Kr09_1})
and from Theorem \ref{theorem 9.26.1}, that the solution
in Theorem \ref{theorem 8.8.3} is independent of $p$
like in Theorem \ref{theorem 9.26.1}
if $\gamma$ is chosen as in Theorem \ref{theorem 9.26.1}
and $\lambda\geq\hbox{\rm RHS of }\eqref{9.27.1}+1$.
 
\end{remark}

\mysection{Differentiating compositions
of generalized functions with differentiable functions}
                                             \label{section 6.9.2}

  Let $\cD$ be the space 
of
generalized functions on   $\bR^{d}$. We need a formula
for $u_{t}(x+x_{t})$ where $u_{t}$ behaves like a function
from $\cW^{1}_{p}$ and $x_{t}$ is  an  $\bR^{d}$-valued
differentiable function. The formula is absolutely natural
and probably well known. We refer the reader to
\cite{Kr09_4} where such a formula is derived in
a much more general setting of stochastic processes. 
 Recall that for any $v\in\cD$
and $\phi\in C^{\infty}_{0}$ the function $(v,\phi(\cdot-x))$
is infinitely differentiable with respect to $x$, so that 
the sup in \eqref{11.16.2} below is measurable.

\begin{definition} 
 Denote by $\frD(S,T)$                   \label{def 10.25.1}
 the set of all $\cD$-valued 
functions $u$ (written 
as $u_{t}(x)$ in a common abuse of notation)
on $[S,T]\cap\bR $ such that, for any $\phi\in C_{0}^{\infty}$,
  the function $(u_{t},\phi)$ is measurable.
Denote by $\mathfrak{D}^{1}(S,T)$ the subset of $\frD(S,T)$
consisting of $u$ such that,  
  for any  
  $\phi\in C_{0}^{\infty}$, $R\in(0,\infty)$,  and finite
$t_{1},t_{2}\in[S,T]$ such that $t_{1}<t_{2}$  we have
\begin{equation}
                                            \label{11.16.2}
\int_{t_{1}}^{t_{2}}\sup_{ |x|\leq R}|(u_{t} ,
\phi(\cdot-x))| \,dt<\infty.
\end{equation}
 
\end{definition}

\begin{definition} 
                                           \label{def 10.25.3}
Let $f,u\in\mathfrak{D}(S,T)$.
 We say that the equation
\begin{equation}
                                           \label{11.16.3}
\partial_{t}u_{t}(x)=f_{t}( x) ,\quad t\in[S,T]\cap\bR,
\end{equation}
holds {\em in the sense of distributions\/} if 
$ f \in\mathfrak{D}^ {1}(S,T)$  and
for
 any $\phi\in C_{0}^{\infty}$  for all finite $s,t\in[S,T]$
we have
$$
(u_{t } ,\phi)=(u_{s} ,\phi)+\int_{s}^{t}
(f_{r},\phi)\,dr .
$$
\end{definition}

Let $x_{t}$ be an $\bR^{d}$-valued function given by
$$
x _{t}=\int_{0}^{t}\hat b _{s}\,ds,
$$
where $\hat b_{s}$ is an $\bR^{d}$-valued locally integrable 
function on $\bR$. Here is the formula.
 
\begin{theorem}
                                   \label{theorem 11.16.5}
Let $f,u\in\mathfrak{D}(S,T) $.
Introduce
$$
v_{t}(x)=u_{t}(x+x_{t})
$$
and assume that \eqref{11.16.3} holds
 (in the sense of distributions). Then
$$
\partial_{t}v_{t}(x)=  f_{t}( x+x_{t})+\hat b^{i}_{t}D_{i}v_{t}(x)
,\quad t\in[S,T]\cap\bR
$$
(in the sense of distributions).
\end{theorem}

\begin{corollary}
                                              \label{corollary 7.3.1}
Under the assumptions of Theorem \ref{theorem 11.16.5}
for any $\eta\in C^{\infty}_{0}$ we have
$$
\partial_{t}[u_{t}(x)\eta(x-x_{t})]=
f_{t}(x)\eta(x-x_{t})-u_{t}(x)
\hat{b}^{i}_{t}D_{i}\eta(x-x_{t}) ,\quad t\in[S,T]\cap\bR.
$$

\end{corollary}

Indeed, what we claim is that for any $\phi\in C^{\infty}_{0}$
and finite $s,t\in[S,T]$
$$
((u_{t } \phi)(\cdot+x_{t }),\eta)=( u_{s} \phi ,\eta)
+\int_{s}^{t} \big(\big[f_{r}\phi 
+ \hat{b}^{i}_{r}D_{i}(u_{r}\phi) \big]
(\cdot+x_{r}),\eta\big)\,dr.
$$
 However, to obtain this result it suffices
to write down an obvious equation for $u_{t}\phi$, then use 
Theorem \ref{theorem 11.16.5} and, finally,
use Definition \ref{def 10.25.3}
to interpret the result.

\mysection{Proof of Theorem \protect\ref{theorem 3.11.1}}
                                           \label{section 6.9.3} 

Throughout this section we suppose that  
Assumptions \ref{assumption 2.7.2},
\ref{assumption 3.11.1}, and \ref {assumption 2.4.01} 
are satisfied (with a $\gamma\in(0,1]$) and
 start with analyzing the  integral in
\eqref{3.16.7}. Recall that $q$ was introduced before
Assumption \ref{assumption 3.11.1}.

\begin{lemma}
                                        \label{lemma 9.1.1}
Let $1\leq r<p$ and
\begin{equation}
                                                    \label{9.1.1}
 \eta:=1+\frac{d}{p}-\frac{d}{r}\geq0
\end{equation}
with strict inequality if $r=1$. Then for any $U\in\cL_{r}$
and $\varepsilon>0$ there exist
$V^{j}\in\cL_{p}$, $j=0,1,...,d$, such that
$U=D_{i}V^{i}+V^{0}$ and
\begin{equation}
                                                    \label{9.1.2} 
\sum_{j=1}^{d}\|V^{j}\|_{\cL_{p}}\leq
N(d,p,r)\varepsilon^{\eta/(1-\eta)}\|U\|_{\cL_{r}},\quad
\|V^{0}\|_{\cL_{p}}\leq
N(d,p,r)\varepsilon^{-1}\|U\|_{\cL_{r}}.
\end{equation}
  In particular, for any $w\in W^{1}_{p'}$
$$
|(U,w)|\leq N(d,p,r)\|U\|_{\cL_{r}} \|w\|_{W^{1}_{p'}}.
$$
\end{lemma}

Proof. If the result is true for $\varepsilon=1$, then
for arbitrary $\varepsilon>0$ it is easily obtained by scaling.
Thus let $\varepsilon=1$
and denote by $R_{0}(x)$ the kernel of $(1-\Delta)^{-1}$.
For $i=1,...,d$  set  $R_{i}=-D_{i}R_{0}$. One knows that
$R_{j}(x)$ decrease exponentially fast as $|x|\to\infty$ and
$$
|R_{j}(x)|\leq\frac{N}{|x|^{d-1}},\quad j=0,1,...,d.
$$
Define
$$
V^{j}=R_{j}*U,\quad j=0,1,...,d.
$$
If $r=1$, one obtains \eqref{9.1.2} from Young's 
inequality since, owing to  the strict
inequality in \eqref{9.1.1} we have $p<d/(d-1)$,
so that $R_{j}\in\cL_{p}$. If $r>1$, then for
 $\nu$ defined by
$$
\frac{1}{p}=\frac{1}{r}-\frac{\nu}{d}
$$
we have $\nu\in(0,1]$, so that
$$
|R_{j}(x)|\leq\frac{N}{|x|^{d-\nu}},\quad j=0,1,...,d,
$$
and we obtain \eqref{9.1.2}  from 
the Sobolev-Hardy-Littlewood inequality. 
After this it only remains to
 notice that in the sense of generalized functions
$$
D_{i}V^{i}+V_{0}=R_{0}*U-\Delta R_{0}*U=U.
$$
The lemma is proved.

Observe that by H\"older's inequality    
for $r=pq/(p+q)$ ($\in[1,p)$ due to $q\geq p'$, see
\eqref{8.11.3}) we have
$$
\|h v\|_{\cL_{r}}\leq  \|h\|_{\cL_{q}}\|v\|_{\cL_{p}}.
$$
Furthermore, if $r=1$, then $q=p'>d$ (see
\eqref{8.11.3}), $p<d/(d-1)$, and $\eta>0$.
In this way
 we come to the following.

\begin{corollary}
                                            \label{corollary 9.1.1}
Let $h\in\cL_{q}$, $v\in
\cL_{p}$,   and $w\in W^{1}_{p'}$. Then  
for any $\varepsilon>0$ there exist
$V^{j}\in\cL_{p}$, $j=0,1,...,d$, such that
$h v=D_{i}V^{i}+V^{0}$ and
$$
\sum_{j=1}^{d}\|V^{j}\|_{\cL_{p}}\leq
N(d,p)\varepsilon^{(q-d)/d}\|h\|_{\cL_{q}}\|v\|_{\cL_{p}}, 
$$
$$
\|V^{0}\|_{\cL_{p}}\leq
N(d,p)\varepsilon^{-1}\|h\|_{\cL_{q}}\|v\|_{\cL_{p}}.
$$
  In particular,
\begin{equation}
                                                         \label{7.1.1}
|(hv,w)|\leq N(d,p)\|h\|_{\cL_{q}}\|v\|_{\cL_{p}}\|w\|_{W^{1}_{p'}}.
\end{equation}
\end{corollary}

\begin{lemma}
                                        \label{lemma 9.1.2}
Let $h\in\cL_{q}$ and
 $u\in W^{1}_{p}$. Then for any $\varepsilon>0$ we have
\begin{equation}
                                                         \label{9.1.5}
\|hu\|_{\cL_{p}}\leq N(d,p)\|h\|_{\cL_{q}}
\big(\varepsilon^{(q-d)/d}
\|Du\|_{\cL_{p}}+\varepsilon^{-1}\|u\|_{\cL_{p}}\big).
\end{equation}
\end{lemma}

Proof. As above it suffices to concentrate on $\varepsilon=1$.
In case $q>p$ 
observe that by H\"older's inequality 
$$
\| hu\|_{\cL_{p}}
\leq\|h\| _{\cL_{q}}\| u\| _{\cL_{s}},
$$
where $s=pq/(q-p)$. After that it only remains to use
embedding theorems (notice that $ 1-d/p\geq -d/s$
since $q\geq d$).
 In the remaining case $q=p$,
which happens only if $p>d$ (see \eqref{8.11.3}). In that case the above
estimate  remains  true if we set $s=\infty$.
The lemma is proved. 

Before we extract some consequences from the lemma we 
take  a  nonnegative
$ \xi\in C^{\infty}_{0}(B_{\rho_{1}})$
  with unit integral and 
define
$$
\bar{b}_{s}(x)=\int_{B_{\rho_{1}}}\xi(y) b_{s}(x-y) \,dy,\quad
\bar{\gb}_{s}(x) =\int_{B_{\rho_{1}}}\xi(y) \gb_{s}(x-y) \,dy,
$$
\begin{equation}
                                                      \label{6.28.3}
\bar{c}_{s}(x)=\int_{B_{\rho_{1}}}\xi(y) c_{s}(x-y) \,dy.
\end{equation}
We may assume that $|\xi|\leq N(d)\rho_{1}^{-d}$.

One obtains the first   assertion  of the following corollary
from \eqref{7.1.1} 
by observing that
$$
\|I_{B_{\rho_{1}}(x_{t})}(b_{t}-\bar{b}_{t}(x_{t}))\|_{\cL_{q}}^{q}
= \int_{B_{\rho_{1}}(x_{t})}| 
b_{ t}-\bar{b}_{t}(x_{t})|^{q}\,dx  
$$
$$
=\int_{B_{\rho_{1}}(x_{t})}\big|\int_{B_{\rho_{1}}(x_{t})}
[b_{t}(x)-b_{t}(y)]\xi(x_{t}-y)\,dy\big|^{q}\,dx
$$
$$
\leq N\int_{B_{\rho_{1}}(x_{t})}\big|\rho_{1}^{-d}
\int_{B_{\rho_{1}}(x_{t})}
|b_{t}(x)-b_{t}(y)|\,dy\big|^{q}\,dx
$$
\begin{equation}
                                                \label{6.11.2}
\leq N\rho_{1}^{-d}\int_{B_{\rho_{1}}(x_{t})} 
\int_{B_{\rho_{1}}(x_{t})}
|b _{t}(x)-b_{t}(y)|^{q}\,dy \,dx\leq 
N\rho_{1}^{-d}KI_{q>d}+N\gamma  .
\end{equation}

The second assertion follows from  
 estimates like \eqref{6.11.2} and \eqref{9.1.5} where
one chooses 
$\varepsilon$ appropriately if $q>d$.

\begin{corollary}
                                  \label{corollary 6.27.2}
Let $u\in\bW^{1}_{p}(S,T)$,
 let $x_{s}$ be
an $\bR^{d}$-valued measurable function, and let
$\eta\in C^{\infty}_{0}(B_{\rho_{1}})$. Set
$\eta_{s}(x)=\eta(x-x_{s})$,
$$
K_{1}=\sup|\eta|+\sup|D\eta|.
$$
 Then   on $(S,T) $

(i) For any $w\in W^{1}_{p'}$ and $v\in \cL_{p}$
$$
 (|b _{s}-\bar{b} _{s}(x_{s})|\eta_{s}
v,|w|) 
\leq N(d,p,K) \|\eta_{s}v\| _{\cL_{p}}
\|w\|_{W^{1}_{p'}} ; 
$$

(ii) We have
$$
\|\eta_{s}
|\gb_{s}-\bar{\gb}_{s}(x_{s})|\,u_{s}\|_{\cL_{p}}
+\|\eta_{s}|c_{s}-\bar{c}_{s}(x_{s})|\,u_{s}\|_{\cL_{p}}
$$
$$
\leq N(d,p)\gamma^{1/q}  
 \|\eta_{s}Du_{s}\| _{\cL_{p}}
+N(d,p ,\gamma,\rho_{1},
K,K_{1})\|I_{B_{\rho_{1}}(x_{s})}u_{s}\| _{\cL_{p}} .
$$

(iii) Almost everywhere on  $(S,T) $ we have
\begin{equation}
                                                \label{7.1.3}
(b^{i}_{s}-\bar{b}^{i}_{s}(x_{s}))\eta_{s} D_{i}u_{s}
=D_{i}V^{i}_{s}+V^{0}_{s} ,
\end{equation}
$$
\sum_{j=1}^{d}\|V^{j}\|_{\cL_{p}}\leq
N(d,p)\gamma^{1/q}\|\eta_{s}Du_{s}\|_{\cL_{p}}, 
$$
\begin{equation}
                                                \label{6.11.3}
 \|V^{0}_{s}\|_{\cL_{p}}
\leq N(d,p,\gamma,\rho_{1},K) \|\eta_{s}Du_{s}\|_{\cL_{p}},
\end{equation}
where $V^{j}_{s} $, $j=0,...,d$, are some measurable 
$\cL_{p}$-valued functions on $(S,T) $.

 \end{corollary}

To prove (iii) observe that one can find a Borel
set $A\subset(S,T)$
of full measure such that  
$I_{A}D_{i}u$, $i=1,...,d$,
 are well defined as $\cL_{p}$-valued Borel measurable
functions. Then \eqref{7.1.3} 
with $I_{A}D_{i}u$ in place of $D_{i}u$
and \eqref{6.11.3} follow  from \eqref{6.11.2},
Corollary \ref{corollary 9.1.1},
 and the fact that the way $V^{j}$ are constructed
uses bounded hence continuous operators and  translates
the measurability of the data into the measurability of the result.
Since we are interested in
\eqref{7.1.3} and \eqref{6.11.3} holding only almost 
everywhere on  $(S,T) $,
 there is no  actual need for the replacement.

\begin{corollary}
                                  \label{corollary 6.27.1}
Let $u\in\bW^{1}_{p}(S,T)$, $R\in(0,\infty)$,
 $\phi\in C^{\infty}_{0}(B_{R})$,
and let finite $S',T'\in(S,T)$ be such that $S'<T'$.
 Then there is a constant $N$ independent of $u$ and $\phi$
such that
\begin{equation}
                                             \label{8.12.3}
\int_{S'}^{T'}
 (|(b^{i}_{s}D_{i}u_{s},\phi)|+
 |(\gb^{i}_{s} u_{s},D_{i}\phi)|+
|(c_{s}u_{s},\phi)|)\,ds\leq N\|u\|_{\bW^{1}_{p}(S,T)}
\|\phi\|_{W^{1}_{p'}},
\end{equation}
so that requirement (i) in Definition \ref{definition 3.20.01}
can be dropped.  

\end{corollary}
 Proof. By having in mind partitions of unity we convince
ourselves that it suffices to prove \eqref{8.12.3}
under the assumption that $\phi$ has support in a ball
$B$ of radius $\rho_{1}$. Let $x_{0}$
be the center of $B$ and set $x_{s}\equiv x_{0}$. 
Observe that  the estimates from Corollary
\ref{corollary 6.27.2} imply that
$$
|(\gb^{i}_{s} u_{s},D_{i}\phi)|\leq
|( \gb^{i}_{s}-\bar{\gb}^{i}_{s}(x_{0}))
 u_{s},D_{i}\phi)|
+|\bar{\gb}^{i}_{s}(x_{0})(  u_{s},D_{i}\phi)| 
$$
$$
\leq N\|u_{s}\|_{W^{1}_{p}}\|\phi\|_{W^{1}_{p'}}+
|\bar{\gb} _{s}(x_{0})|\,\|u_{s}\|_{W^{1}_{p}}\|\phi\|_{W^{1}_{p'}}.
$$
By recalling
Assumption
\ref{assumption 2.7.2} (iii) and H\"older's inequality we get
$$
\int_{S'}^{T'}|(\gb^{i}_{s} u_{s},D_{i}\phi)|\,ds
\leq N\|u\|_{\bW^{1}_{p}(S,T)}
\|\phi\|_{W^{1}_{p'}}.
$$

Similarly the integrals of $ |(b^{i}_{s}D_{i}u_{s},\phi)|$
and
$|(c_{s}u_{s},\phi)|$ are estimated and the corollary is proved.

Since bounded linear operators
are continuous we obtain the following.
 \begin{corollary}
                                    \label{corollary 3.23.1}
Let $\phi\in C^{\infty}_{0}$, $T\in(0,\infty)$. Then
 the operators
$$
u_{\cdot}\to \int_{0}^{\cdot}
(b^{i}_{t}D_{i}u_{t},\phi)\,dt,\quad
u_{\cdot}\to \int_{0}^{\cdot}
(\gb^{i}_{t}u_{t},D_{i}\phi)\,dt, \quad
u_{\cdot}\to \int_{0}^{\cdot}
(c_{t}u_{t}, \phi)\,dt
$$
are continuous as operators from $\bW^{1}_{p}(\infty)$ to
$\cL_{p}([-T ,T ])$.
\end{corollary}  
This result will be used in Section
\ref{section 6.9.5}.

Before we continue with the proof of Theorem \ref{theorem 3.11.1},
 we notice that, if $u \in 
\cW^{ 1}_{p }(S,T)$, then as we know (see, for instance,
 Theorem 2.1 of \cite{Kr09_3}),
the function $u_{t}$
is a continuous $\cL_{p}$-valued function on
$[S,T]\cap\bR$.  

Now we are ready to prove Theorem \ref{theorem 3.11.1}
in a particular case.

\begin{lemma}
                                        \label{lemma 3.11.1}
Let  $\gb^{i}$, $b^{i}$, and $c$ 
be independent of $x$ and let $S=-\infty$. 
Then the assertion of Theorem \ref{theorem 3.11.1}
holds, naturally, with $\lambda_{0}=\lambda_{0}(d,\delta, p,
\rho_{0} )$ (independent of $\rho_{1}$ and $K$).

\end{lemma}

Proof.  First let  $c\equiv0$.
We want to use Theorem \ref{theorem 11.16.5} to get rid
of the first order terms. Observe that \eqref{2.6.4} reads as
\begin{equation}
                                                \label{6.28.1}
\partial_{t}u_{t}= D_{i}(a^{ij}_{t}D_{j}u_{t}+[\gb^{i}_{t} + 
b^{i}_{t}]u_{t}+f^{i}_{t})+f^{0}_{t}- \lambda u_{t},
\quad t\leq T.
\end{equation}
 
Recall that from the start
(see Definition \ref{definition 3.20.01})
 it is assumed that $u \in  
\bW^{ 1}_{p}(T)$. Then
one can find a Borel set $A\subset(-\infty,T)$
of full measure such that $I_{A}f^{j}$, $j=0,1,...,d$, and
$I_{A}D_{i}u$, $i=1,...,d$,
 are well defined as $\cL_{p}$-valued Borel
functions satisfying
$$
\int_{-\infty}^{T}I_{A}\big(\sum_{j=0}^{d}\|f^{j}_{t}\|^{p}_{\cL_{p}}
+ \|Du_{t}\|^{p}_{\cL_{p}}\big)\,dt<\infty.
$$
Replacing $f^{j}$ and $D_{i}u$ in \eqref{6.28.1} with
$I_{A}f^{j}$ and $I_{A}D_{i}u$, respectively, will not affect
\eqref{6.28.1}. Similarly one can treat the term 
$h_{t}= (\gb^{i}_{t}+b^{i}_{t})u_{t}$
for which
$$
\int_{S'}^{T'}\|h_{t}\|_{\cL_{p}}\,dt<\infty
$$
for each finite $S',T'\in(-\infty ,T]$, owing to Assumption 
\ref{assumption 2.7.2} and 
the fact that $u\in\bL_{p}( T)$.

After these replacements all terms on the right in \eqref{6.28.1}
will be of   class $\frD^{1}(-\infty,T)$   since  
 $a$  is bounded. 
This allows us to apply Theorem \ref{theorem 11.16.5}
and for
$$
B_{t}^{i}=\int_{0}^{t}
(\gb^{i}_{s}+b^{i}_{s})\,ds,\quad 
 \hat{u}_{t}(x)=u_{t}(x-B_{t})
$$
obtain that
\begin{equation}
                                               \label{6.28.8}
\partial_{t}\hat{u}_{t}=D_{i}(\hat{a}^{ij}_{t}D_{j}\hat{u}_{t} )
 - \lambda \hat{u}_{t}+D_{i}\hat{f}^{i}_{t}+\hat{f}^{0}_{t},
\end{equation}
where
$$
(\hat{a}^{ij}_{t},
\hat{f}^{j}_{t})(x)=(a^{ij}_{t}, f^{j}_{t})(x-B_{t}).
$$

Obviously, $\hat{u}$ is in $\bW^{1}_{p}(T)$ and its norm
coincides with that of $u$. 
Equation \eqref{6.28.8}  
shows that $\hat{u}\in\cW^{1}_{p}(T)$.

By Theorem 4.4 and Remark 2.4 of \cite{Kr07} there exist
$\gamma=\gamma(d,\delta, p)$ and $\lambda_{0}=
\lambda_{0}(d,\delta, p,
\rho_{0} )$ such that
if $\lambda\geq\lambda_{0}$, then
\begin{equation}
                                                      \label{8.4.1}
\|D\hat{u}\|_{\bL_{p}(T)} +\lambda^{1/2}\|
\hat{u}\|_{\bL_{p}(T)}\leq N\big(\sum_{i=1}^{d}
\|\hat f^{i}\|_{\bL_{p}(T)}
+\lambda^{-1/2}\|\hat f^{0}\|_{\bL_{p}(T)}
\big).
\end{equation}
Actually,  Theorem 4.4 of \cite{Kr07} is proved there
only for $T=\infty$, but it is a standard fact that such an estimate
implies what we need for any $T$ (cf. the proof of Theorem
6.4.1 of \cite{Kr 2008}).
Since the norms in $\cL_{p} $ and $W^{1}_{p}$
are translation invariant, \eqref{8.4.1}
implies \eqref{3.11.2}  and finishes the proof of the lemma
in case $c\equiv0$.

Our next step is to abandon the condition $c\equiv0$ but assume that
for an
$S>-\infty$ we have $u_{t}=f^{j}_{t}=0$ for $t\leq S$. 
Observe that without loss of generality we may assume that $T
<\infty$.
In that case introduce
$$
\xi_{t}=\exp(\int_{S}^{t}c_{s}\,ds).
$$
Then    we have $v:=
\xi u\in \bW^{1}_{p}(T)$ and
$$
\partial_{t}v_{t}= D_{i}(a^{ij}_{t}D_{j}v_{t}+[\gb^{i}_{t} + 
b^{i}_{t}]v_{t}+\xi_{t}f^{i}_{t})+\xi_{t}f^{0}_{t}- \lambda v_{t},
\quad t\leq T.  
$$
By the above result for all   $T'\leq T$ 
$$
\int_{-\infty}^{T'}\xi_{t}^{p}\|Du_{t}\|_{\cL_{p}}^{p}\,dt
+\lambda^{p/2}
\int_{-\infty}^{T'}\xi_{t}^{p}\|u_{t}\|_{\cL_{p}}^{p}\,dt
$$
\begin{equation}
                                                     \label{8.8.1}
\leq N_{1}\sum_{i=0}^{d}
\int_{-\infty}^{T'}\xi_{t}^{p}\|f^{i}_{t}\|_{\cL_{p}}^{p}\,dt
+N_{1}\lambda^{-p/2}
\int_{-\infty}^{T'}\xi_{t}^{p}\|f^{0}_{t}\|_{\cL_{p}}^{p}\,dt.
\end{equation}
We multiply both part of \eqref{8.8.1} by $pc_{T'}\xi_{T'}^{-p}$
and integrate with respect to $T'$ over $(S,T)$. We use integration 
by parts observing that both parts vanish at $T'=S$. Then we obtain
$$
\int_{-\infty}^{T}\|Du_{t}\|_{\cL_{p}}^{p}\,dt
+\lambda^{p/2}
\int_{-\infty}^{T}\|u_{t}\|_{\cL_{p}}^{p}\,dt
$$
$$
-\xi_{T}^{-p}\int_{-\infty}^{T}\xi_{t}^{p}\|Du_{t}\|_{\cL_{p}}^{p}\,dt
-\xi_{T}^{-p}\lambda^{p/2}
\int_{-\infty}^{T}\xi_{t}^{p}\|u_{t}\|_{\cL_{p}}^{p}\,dt
$$
$$
\leq N_{1}\sum_{i=0}^{d}
\int_{-\infty}^{T}\|f^{i}_{t}\|_{\cL_{p}}^{p}\,dt
+N_{1}\lambda^{-p/2}
\int_{-\infty}^{T}\|f^{0}_{t}\|_{\cL_{p}}^{p}\,dt.
$$
$$
-\xi_{T}^{-p} N_{1}\sum_{i=0}^{d}
\int_{-\infty}^{T}\xi_{t}^{p}\|f^{i}_{t}\|_{\cL_{p}}^{p}\,dt-
\xi_{T}^{-p}N_{1}\lambda^{-p/2}
\int_{-\infty}^{T}\xi_{t}^{p}\|f^{0}_{t}\|_{\cL_{p}}^{p}\,dt.
$$
By adding up this inequality with \eqref{8.8.1} with $T'=T$
multiplied by $\xi_{T}^{-p}$ we obtain \eqref{3.11.2}.

The last step is
 to avoid assuming that $u_{t}=0$ for large negative $t$.
In that case we find a sequence $S_{n}\to-\infty$
such that $u_{S_{n}}\to0$  in $W^{1}_{p}$ and denote by $v^{n}_{t}$
the unique solution of class $\bW^{1}_{p}((0,1)\times\bR^{d})$
of the heat equation $\partial v^{n}_{t}=\Delta v^{n}_{t}$
with initial condition $u_{S_{n}}$. After that we modify $u_{t}$
and the coefficients of $L_{t}$ for $t\leq S_{n}$ as in the proof
of Theorem \ref{theorem 8.8.1} by taking there $v^{n}_{t}$ and 
$S_{n}$ in place of $v_{t}$ and $S$, respectively. Then by the above result
  we obtain
$$
\lambda\| u\|^{2}_{\bL _{p}(S_{n},T)}+
\|D u\|^{2}_{\bL _{p}(S_{n},T)}
\leq N\big(\sum_{i=1}^{d}\|\tilde f^{i}\|^{2}_{\bL _{p}(T)}
+\lambda^{-1}\|\tilde f^{0}\|^{2}_{\bL _{p}(T)} \big),
$$
$$
\leq N\big(\sum_{i=1}^{d}\| f^{i}\|^{2}_{\bL _{p}(T)}
+\lambda^{-1}\|  f^{0}\|^{2}_{\bL _{p}(T)} \big)
+N(1+\lambda^{-1})\|u_{S_{n}}\|_{W^{1}_{p}}^{p}.
$$
By letting $n\to\infty$ we come to \eqref{3.11.2} and
the lemma is proved.
\begin{remark}
In \cite{Kr07} the assumption corresponding
to Assumption \ref{assumption 2.4.01} is much
weaker since in the corresponding
counterpart of \eqref{2.9.5}
there is no supremum over $x\in\bR^{d}$.
We need  our stronger assumption because
we need $a^{ij}_{t}(x-B_{t})$ to satisfy
the assumption in \cite{Kr07} for any function $B_{t}$.

\end{remark}

To proceed further we need a construction. Recall that $\bar{\gb}$
and
$\bar{b}$ are introduced in \eqref{6.28.3}.
 From Lemma 4.2 of \cite{Kr09_1} and Assumption
\ref{assumption 3.11.1} it follows  that, for 
$h_{ t}=\bar{\gb}_{ t},\bar{b}_{t}$, it holds that  
 $|D^{n}h_{ t}|\leq \kappa_{n} $, where
$\kappa_{n}=\kappa_{n}
(n, d,p,\rho_{1},K)\geq1$ and $ D^{n}h_{ t }$
is any derivative of $h_{ t}$ of order $n\geq1$
with respect to $x$. By Corollary 4.3 of \cite{Kr09_1}
we have $|h_{ t}( x )|\leq K(t)(1+|x|)$,
where  the function $K(t) $
is locally 
  integrable with respect to $t$ on $\bR$.
Owing to these properties,  for any  
$ (t_{0},x_{0}) \in \bR^{d+1 }$,
the equation
$$
x_{t}=x_{0}-\int_{t_{0}}^{t}(\bar{\gb}_{ s}+\bar{b}_{  s})
( x_{s})\,ds,\quad t \geq t_{0} ,
$$
  has a unique solution 
$x_{t}=x_{t_{0},x_{0},t} $. 

Next, for $i=1,2$ set $\chi^{(i)}(x)$ to be the indicator function
of $B_{\rho_{1}/i}$ and introduce
$$
\chi^{(i)}_{t_{0},x_{0},t}(x)=\chi^{(i)}(x-x_{t_{0},x_{0},t}).
$$

Here is a crucial estimate.
\begin{lemma}
                                     \label{lemma 3.14.1}
Suppose that  
Assumptions \ref{assumption 2.7.2},
\ref{assumption 3.11.1}, and \ref {assumption 2.4.01} 
are satisfied with a $\gamma\in(0,\gamma(d,p,\delta)]$,
where $\gamma(d,p,\delta)$ is taken from Lemma \ref{lemma 3.11.1}.
Take $ (t_{0},x_{0}) \in 
\bR^{d+1}$ and assume that $t_{0}<T$ and that we are given a function
$u  $ which is a solution of \eqref{2.6.4} with $S=t_{0}$,
with  zero initial condition,  
   some
$f^{j} \in\bL_{p}(t_{0},T)$, and $\lambda\geq\lambda_{0}$,
 where $\lambda_{0}=\lambda_{0}(d,\delta, p, 
\rho_{0} )$ is taken from Lemma~\ref{lemma 3.11.1}.
 Then  
$$
 \lambda \|\chi^{(2)}_{t_{0},x_{0}}u \|^{2}_{\bL _{p}(t_{0},T)}+
\|\chi^{(2)}_{t_{0},x_{0}}Du\|^{2}_{\bL _{p}(t_{0},T)} 
$$
$$
\leq N\sum_{i=1}^{d}\|\chi^{(1)}_{t_{0},x_{0}}
f^{i}\|^{2}_{\bL _{p}(t_{0},T)}
+N\lambda^{-1 }\|\chi^{(1)}_{t_{0},x_{0}}f^{0}\|^{2}_{\bL _{p}(t_{0},T)}
$$
$$
+N\gamma ^{2/q} \| 
\chi^{(1)}_{t_{0},x_{0}} Du
\|_{\bL_{p}(t_{0},T)}^{2}+
 N^{*} \lambda^{-1}\| 
\chi^{(1)}_{t_{0},x_{0}} Du
\|_{\bL_{p}(t_{0},T)}^{2}
$$
\begin{equation}
                                       \label{3.14.2}
+ N^{*}  \|
\chi^{(1)}_{t_{0},x_{0}}  u
\|_{\bL_{p}(t_{0},T)}^{2}
+N^{*}\lambda^{-1}\sum_{i=1}^{d}\|\chi^{(1)}_{t_{0},x_{0}}
f^{i}\|^{2}_{\bL _{p}(t_{0},T)},
\end{equation}
where and below  in the proof
 by $N$ we denote generic constants depending only
on $d,\delta$,    and $p$ and by $N^{*}$ constants depending only
on the same objects, $\gamma$, $\rho_{1}$,  and $K$.
\end{lemma}

Proof. Shifting the origin  
allows us to assume that $t_{0}=0$ and $x_{0}=0$.
  With this stipulations we will drop the subscripts $t_{0},
x_{0}$. 

Fix  a 
$ \zeta\in C^{\infty}_{0} $ with support in $B_{\rho_{1}}$
and such that $\zeta =1$ on $B_{\rho_{1}/2}$ and $0\leq\zeta
\leq1$. 
Set $x_{t}=x_{0,0,t}$,
$$
  \hat{\gb}_{ t} =\bar{\gb}_{ t}( x_{ t}) ,\quad
\hat{b}_{ t} =\bar{b}_{  t}( x_{ t}),\quad
\hat{c}_{ t} =\bar{c}_{  t}( x_{ t})
$$
$$
\eta_{ t}( x)=\zeta(x-x_{ t} ), 
\quad v_{ t}( x)= u_{t}( x)  \eta_{ t}( x).
$$
The most important property of $\eta_{t}$ is that
$$
\partial_{t}
\eta_{t}=(\hat{\gb}^{i}_{t}+\hat{b}^{i}_{t})D_{i}\eta_{t}.
$$
Also
observe  for the later that we may assume that
\begin{equation}
                                                   \label{6.16.1}
\chi^{(2)}_{t}\leq\eta_{ t}\leq \chi^{(1)}_{t},\quad
|D\eta_{ t}|\leq N\rho_{1}^{-1 }\chi^{(1)}_{t},
\end{equation}
where $\chi^{(i)}_{t}=\chi^{(i)}_{0,0,t}$ and $N=N(d)$.

By Corollary \ref{corollary 7.3.1} (also see the argument before
\eqref{6.28.8})
we obtain
that  for finite $t\in[0,T]$ 
$$
\partial_{t}
v_{ t} =D_{i}(\eta_{ t}a^{ij}_{t}D_{j}u_{t}+\gb^{i}_{t}v_{ t})
-(a^{ij}_{t}D_{j}u_{t}+\gb^{i}_{t}u_{t})
D_{i}\eta_{ t}
$$
$$
+b^{i}_{t}\eta_{ t}D_{i} u_{t} 
-(c_{t}+\lambda)  v_{ t}
+D_{i}(f^{i}_{t}\eta_{ t})-f^{i}_{t}D_{i}\eta_{ t}
+f^{0}_{t}\eta_{ t}
+
(\hat{\gb}^{i}_{ t} +\hat{b}^{i}_{ t} )u_{t}
D_{i} \eta_{ t}.
$$
We transform  this further by noticing that  
$$
\eta_{ t}a^{ij}_{t}D_{j}u_{t}= a^{ij}_{t}D_{j}v_{ t}-
a^{ij}_{t}u_{t}D_{j}\eta_{ t}.
$$

To deal with the term $b^{i}_{t}\eta_{ t}D_{i} u_{t} $
we use Corollary \ref{corollary 6.27.2} and find 
the corresponding functions $V^{j}_{t}$.
Then simple arithmetics show that
$$
\partial_{t}v_{ t}=
 D_{i}\big(a^{ij}_{t}D_{j}v_{ t}+
\hat{\gb}^{i}_{ t}v_{ t}\big)- (\hat{c}_{t}+\lambda) v_{ t}
+\hat{b}^{i}_{ t}
D_{i}v_{ t}+D_{i}\hat{f}^{i}_{ t}
+\hat{f}^{0}_{ t},
$$  
where
$$
\hat{f}^{0}_{ t}=f^{0}_{t}\eta_{ t}-f^{i}_{t}D_{i}\eta_{ t}
 -a^{ij}_{t}(D_{j}u_{t})D_{i}\eta_{ t}
+( \hat{\gb}^{i}_{ t}-\gb^{i}_{t})
u_{t}D_{i} \eta_{ t} 
+V^{0}_{t}+(\hat{c}_{t}-c_{t})u_{t}\eta_{t},
$$ 
$$
\hat{f}^{i}_{ t}=f^{i}_{t}\eta_{ t}-
a^{ij}_{t}u_{t}D_{j}\eta_{ t}+(\gb^{i}_{t}-
\hat{\gb}^{i}_{ t})
u_{ t}\eta_{ t}+V^{i}_{t},\quad i=1 ,..,d.
$$

It we extend $u_{t}$ and $f^{j}_t$ as zero
for $t<0$, then it will be seen from
  Lemma~\ref{lemma 3.11.1} that for $\lambda\geq\lambda_{0}$
\begin{equation}
                                                   \label{3.13.1}
 \lambda\|v  \|^{2}_{\bL_{p}(0,T)}+
\|Dv \|^{2}_{\bL_{p}(0,T)}\leq 
N
\sum_{i=1}^{d}\|\hat{f}^{i} \|^{2}_{\bL_{p}(0,T)}
+ N\lambda^{-1}\|\hat{f}^{0} \|^{2}_{\bL_{p}(0,T)}.
\end{equation}
Recall that here
  and below by $N$ we denote generic constants
depending only on $d,\delta$, and $p$.

Now we start estimating the right-hand side of \eqref{3.13.1}.  
First we deal with $\hat{f}^{i}_{ t}$.
Recall \eqref{6.16.1} and use Corollary \ref{corollary 6.27.2}
to get  
\begin{equation}
                                          \label{3.14.1}
 \| (\gb^{i}_{t}-
\hat{\gb}^{i}_{ t})
u_{ t}\eta_{ t}\|^{2}_{\cL_{p}}
\leq N \gamma^{2/q} \|\chi^{(1)}_{t}
 Du_{t}  \|^{2}_{\cL_{p}}+
 N^{*}\|\chi^{(1)}_{t}
  u_{t}\|^{2}_{\cL_{p}} 
\end{equation}  
(we remind the reader that
 by $N^{*}$ we denote generic constants depending
only on $d$, $\delta$, $p$,  $\gamma$, $\rho_{1}$,  and $K$).
By adding that
$$
\|a^{ij} u D_{j}\eta \|^{2}
_{\bL_{p}(0,T)}\leq N^{*}
\|\chi^{(1)}_{\cdot}u \|^{2}_{\bL_{p}(0,T)},
$$
we derive from \eqref{6.11.3}  and \eqref{3.14.1}
that  
$$
\sum_{i=1}^{d}\|\hat{f}^{i} \|^{2}_{\bL_{p}(0,T)}
\leq N\sum_{i=1}^{d}\|\chi^{(1)}_{\cdot}f^{i}  
\|^{2}_{\bL_{p}(0,T)}
$$
\begin{equation}
                                          \label{3.14.4}
+N\gamma ^{2/q} \| 
\chi^{(1)}_{\cdot}  Du
\|_{\bL_{p}(0,T)}^{2}+ N^{*}\|
\chi^{(1)}_{\cdot}  u
\|_{\bL_{p}(0,T)}^{2}.
\end{equation}

While estimating $\hat{f}^{0}$ we 
use \eqref{6.11.3} again and
observe that  we can deal with $( \hat{\gb}^{i}_{ t}-\gb^{i}_{t})
u_{t}D_{i} \eta_{ t}$ and $( c_{t}-
\hat{c}_{ t})
u_{ t}\eta_{ t}$ as in \eqref{3.14.1} this time without paying
too much attention to the dependence of our constants on $
\gamma$, $\rho_{1}$, and $K$
and obtain that
$$
\|( \hat{\gb}^{i} -\gb^{i} )
u D_{i} \eta \|_{\bL_{p}(0,T)}^{2}+
  \|( c_{t}-
\hat{c}_{ t})
u_{ t}\eta_{ t}\|^{2}_{\cL_{p}}
$$
$$
\leq  N^{*}(\|\chi^{(1)}_{\cdot} Du
\|_{\bL_{p}(0,T)}^{2}+\|\chi^{(1)}_{\cdot}u\|_{\bL_{p}(0,T)}^{2}).
$$
By estimating also roughly the remaining terms in $\hat{f}^{0}$
and combining this with \eqref{3.14.4} and \eqref{3.13.1},
we see that the left-hand side of \eqref{3.13.1}
is less than the right-hand side of \eqref{3.14.2}.
  However,
$$
|\chi^{(2)}_{t}Du_{t}|\leq|\eta_{t}Du_{t}|\leq
|Dv_{t}|+|u_{t}D\eta_{t}|\leq|Dv_{t}|+
N\rho_{1}^{-1}|u_{t} \chi^{(1)}_{t}|
$$
which easily leads to \eqref{3.14.2}.
The lemma is proved.

Next, from  the  result giving ``local" in space estimates
we derive global in space estimates but for functions
having, roughly
speaking, small ``past" support in the time variable.  
In the following lemma $\kappa_{1}$ is the number
introduced before Lemma \ref{lemma 3.14.1}.

\begin{lemma}
                                                \label{lemma 3.14.3}

Suppose that  
Assumptions \ref{assumption 2.7.2},
\ref{assumption 3.11.1}, and \ref {assumption 2.4.01} 
are satisfied with a $\gamma\in(0,\gamma(d,p,\delta)]$,
where $\gamma(d,p,\delta)$ is taken from Lemma \ref{lemma 3.11.1}.
Assume that
$u  $ is a solution of \eqref{2.6.4}
  with $S=-\infty$, some
$f^{j} \in\bL_{p}(T)$, and $\lambda\geq\lambda_{0}$,
 where $\lambda_{0}=\lambda_{0}(d,\delta, p, 
\rho_{0} )$ is taken from Lemma~\ref{lemma 3.11.1}.
 Take a finite $t_{0}\leq T$ 
 and assume that $u_{t}=0$ if $t\leq t_{0}$.
 Then for   $I_{t_{0}}:=
I_{(t_{0},T')}$, where $T'=(t_{0}+\kappa_{1}^{-1})\wedge T$,
we have
$$
\lambda^{p/2}\| I_{t_{0}}u\|^{p}_{\bL _{p} }+
\|I_{t_{0}} Du\|^{p}_{\bL _{p} }
\leq N \sum_{i=1}^{d}\| 
I_{t_{0}}f^{i}\|^{p}_{\bL _{p} }
+N\lambda^{-p/2 }\|I_{t_{0}} f^{0}\|^{p}_{\bL _{p} }
$$
$$
 +N\gamma ^{p/q} \| I_{t_{0}} Du\|_{\bL_{p} }^{p}
+
 N^{*} \lambda^{-p/2}\| I_{t_{0}}  Du\|_{\bL_{p} }^{p}
$$
\begin{equation}
                                       \label{3.14.5} 
+ N^{*} \| I_{t_{0}} u\|_{\bL_{p} }^{p}
+N^{*}\lambda^{-p/2}\sum_{i=1}^{d}\| 
I_{t_{0}}f^{i}\|^{p}_{\bL _{p} },
\end{equation}
where and below in the proof by $N$ we denote generic constants depending only
on $d,\delta$, and $p$ and by $N^{*}$ constants depending only
on the same objects, $\gamma$, $\rho_{1}$, and $K$.

\end{lemma}

Proof. Take $x_{0}\in\bR^{d}$ and use the notation introduced before
Lemma \ref{lemma 3.14.1}. By this lemma with $T'$ in place of $T$
we have 
$$
 \lambda^{p/2} \|I_{t_{0}}
\chi^{(2)}_{t_{0},x_{0}}u \|_{\bL _{p} }^{p}+
\|I_{t_{0}}\chi^{(2)}_{t_{0},x_{0}}Du\|^{p}_{\bL _{p} } 
$$
$$
\leq N\sum_{i=1}^{d}\|I_{t_{0}}\chi^{(1)}_{t_{0},x_{0}}
f^{i}\|^{p}_{\bL _{p} }
+N\lambda^{-p/2 }\|I_{t_{0}}\chi^{(1)}_{t_{0},x_{0}}f^{0}\|^{p}_{\bL _{p} }
$$
$$
+N\gamma ^{p/q} \| I_{t_{0}}
\chi^{(1)}_{t_{0},x_{0}} Du
\|_{\bL_{p} }^{p}+
 N^{*} \lambda^{-p/2}\| I_{t_{0}}
\chi^{(1)}_{t_{0},x_{0}} Du
\|_{\bL_{p} }^{p}
$$
\begin{equation}
                                       \label{8.5.1}
+ N^{*}  \|I_{t_{0}}
\chi^{(1)}_{t_{0},x_{0}}  u
\|_{\bL_{p} }^{p}
+N^{*}\lambda^{-p/2}\sum_{i=1}^{d}\|I_{t_{0}}\chi^{(1)}_{t_{0},x_{0}}
f^{i}\|^{p}_{\bL _{p}}.
\end{equation}

One knows that for each   $t\geq t_{0}$,
the mapping $x_{0}\to x_{t_{0},x_{0},t} $ is a diffeomorphism
with Jacobian determinant given by
$$
\bigg|\frac{\partial x_{t_{0},x_{0},t} }{
\partial x_{0}}\bigg| =\exp\big(-\int_{t_{0}}^{t}\sum_{i=1}^{d} D_{i}
[\bar{\gb} _{  s}^{i}+\bar{ b} _{  s}^{i}]
( x_{t_{0},x_{0},s}) \,ds\big).
$$
By the  way the constant $\kappa_{1}$ is introduced, we have
$$
e^{-N\kappa_{1}(t-t_{0})}\leq \bigg|\frac{\partial x_{t_{0},x_{0},t}}{
\partial x_{0}}\bigg|  \leq e^{N\kappa_{1}(t-t_{0})},
$$
where   $N$ depends only on $d$.
  Therefore, for any nonnegative
Lebesgue measurable function $w(x)$ it holds that
$$
e^{-N\kappa_{1}(t-t_{0})}
\int_{\bR^{d}}w(y)\,dy\leq
\int_{\bR^{d}}w(x_{t_{0},x_{0},t})\,dx_{0}\leq 
e^{N\kappa_{1}(t-t_{0})}\int_{\bR^{d}}w(y)\,dy .
$$
In particular, since
$$
\int_{\bR^{d}}|\chi^{(i)}_{t_{0},x_{0},t}( x)|^{p}\,dx_{0}=
\int_{\bR^{d}}|\chi^{(i)}(x-x_{t_{0},x_{0},t})|^{p}\,dx_{0} ,
$$
we have
$$
e^{-N\kappa_{1}(t-t_{0})}=N^{*}_{i} e^{-N\kappa_{1}(t-t_{0})}
 \int_{\bR^{d}}|\chi^{(i)}(x-y)|^{p}\,dy 
$$
$$
 \leq N^{*}_{i} 
\int_{\bR^{d}}|\chi^{(i)}_{t_{0},x_{0},t}( x)|^{p}\,dx_{0}
\leq N^{*}_{i} e^{N\kappa_{1}(t-t_{0})}
\int_{\bR^{d}}|\chi^{(i)}(x-y)|^{p}\,dy=e^{N\kappa_{1}(t-t_{0})}  ,
$$
where $N^{*}_{i}=|B_{1}|^{-1}
\rho_{1}^{-d}i^{d}$ and $|B_{1}|$ is the volume of  
$B_{1}$. It follows that
$$
\int_{\bR^{d}}|\chi^{(1)}_{t_{0},x_{0},t}( x)|^{p}\,dx_{0}
\leq (N^{*}_{1})^{-1}e^{N\kappa_{1}(t-t_{0})},
$$
$$
(N^{*}_{2})^{-1}e^{-N\kappa_{1}(t-t_{0})}\leq
\int_{\bR^{d}}|\chi^{(2)}_{t_{0},x_{0},t}( x)|^{p}\,dx_{0}.
$$
Furthermore, since $u_{t}=0$ if $ t\leq t_{0}$
and $T'\leq t_{0}+\kappa_{1}^{-1}$,
in evaluating the norms in \eqref{8.5.1} we need not
integrate with respect to $t$ such that $\kappa_{1}(t-t_{0})\geq
1$ or $\kappa_{1}(t-t_{0})\leq0$, so that for all $t$ really involved we have
$$
\int_{\bR^{d}}|\chi^{(1)}_{t_{0},x_{0},t}( x)|^{2}\,dx_{0}
\leq (N^{*}_{1})^{-1}e^{N },\quad
(N^{*}_{2})^{-1}e^{-N }\leq
\int_{\bR^{d}}|\chi^{(2)}_{t_{0},x_{0},t}( x)|^{2}\,dx_{0}.
$$
After this observation it only remains to integrate 
\eqref{8.5.1} through with respect to $x_{0}$ and use the 
fact that $N^{*}_{1}=2^{-d}N^{*}_{2}$.
The lemma is proved.
 
{\bf Proof of Theorem \ref{theorem 3.11.1}}. 
Obviously we may assume that $S=-\infty$. Then
first we show how to choose an appropriate
$\gamma =\gamma (d,\delta,p)\in(0,1]$. 
For one, we take it smaller than the one from Lemma
\ref{lemma 3.11.1}.
Then call $N_{0}$
the constant factor of $\gamma ^{p/q} \| 
I_{t_{0}} Du\|_{\bL_{p}}^{p}$
in \eqref{3.14.5}. We know that $N_{0}=N_{0}(d,\delta,p)$ and we choose
$\gamma \in(0,1]$ so that $N_{0}\gamma ^{p/q}\leq1/2$.
Then under the conditions of Lemma \ref{lemma 3.14.3}  
we have 
$$
\lambda^{p/2}\|I_{t_{0}} u \|^{p}_{\bL _{p}}+
\|I_{t_{0}} Du\|^{p}_{\bL _{p}}
\leq N\sum_{i=1}^{d}\| 
I_{t_{0}}f^{i}\|^{p}_{\bL _{p}} 
+N\lambda^{-p/2 }\|I_{t_{0}} f^{0}\|^{p}_{\bL _{p} } 
$$
\begin{equation}
                                       \label{3.14.6}
+
 N^{*} \lambda^{-p/2}\|  I_{t_{0}} Du\|_{\bL_{p}}^{p} 
+ N^{*}\| I_{t_{0}} u\|_{\bL_{p}}^{p} 
+N^{*}\lambda^{-p/2}\sum_{i=1}^{d}\| 
I_{t_{0}}f^{i}\|^{p}_{\bL _{p}}.
\end{equation}
After $\gamma $ has been fixed we recall that $\kappa_{1}=\kappa_{1}
(d, p,\rho_{1},K)$ and  take a $\zeta\in C^{\infty}_{0}(\bR)$
with support in $(0,\kappa_{1}^{-1})$
such that
\begin{equation}
                                       \label{3.15.1} 
\int_{-\infty}^{\infty}\zeta^{p}(t)\,dt=1.
\end{equation}
 For $s\in\bR$ define
$\zeta^{s}_{t}=\zeta(t-s)$, 
$u^{s}_{t}( x)=u_{t}(x)\zeta^{s}_{t}$. 
Obviously $u^{s}_{t}=0$ if $t\leq s\wedge T$.
Therefore, we can apply 
\eqref{3.14.6}  to $u^{s}_{t}$ with $t_{0}=s\wedge T$
 observing that
$$
\partial_{t}u^{s}_{t}= D_{i}(a^{ij}_{t}D_{j}u^{s}_{t}
+\gb^{i}_{t}u^{s}_{t})+b^{i}_{t}u^{s}_{t}-(c_{t}+\lambda)
u^{s}_{t}+D_{i}(\zeta^{s}_{t}f^{i}_{t})+ 
\zeta^{s}_{t}f^{0}_{t}+(\zeta^{s}_{t})'u_{t} .
$$
Then from \eqref{3.14.6} for $\lambda\geq\lambda_{0}  $, 
where $\lambda_{0}=\lambda_{0}(d,\delta, p, 
\rho_{0} )$ is taken from Lemma~\ref{lemma 3.11.1},  we obtain
$$
 \lambda^{p/2} \|I_{s\wedge T }\zeta^{s}  u \|^{p}_{\bL _{p}}+
\|I_{s\wedge T }\zeta^{s}  Du\|^{p}_{\bL _{p}}
\leq N\sum_{i=1}^{d}\| I_{s\wedge T }
\zeta^{s} f^{i}\|^{p}_{\bL _{p}}
$$
$$
+N\lambda^{-p/2 } \|I_{s\wedge T }\zeta^{s}f^{0}\|^{p}_{\bL _{p}}
+N\|I_{s\wedge T }(\zeta^{s} )'u \|^{p}_{\bL _{p}} 
$$
\begin{equation}
                                       \label{3.14.7} 
+ N^{*} \lambda^{-p/2}\|I_{s\wedge T }\zeta^{s}Du\|_{\bL_{p}}^{p}
+ N^{*}\| I_{s\wedge T }\zeta^{s}  u\|_{\bL_{p}}^{p}
+N^{*}\lambda^{-p/2}\sum_{i=1}^{d}\| 
I_{s\wedge T }\zeta^{s} f^{i}\|^{p}_{\bL _{p}}.
\end{equation}
We 
 integrate through \eqref{3.14.7} with respect to $s\in\bR$, observe
that
$$
I_{s\wedge T<t<[(s\wedge T)+\kappa_{1}^{-1}]
\wedge T}=I_{t< T}
I_{s\wedge T<t<(s\wedge T)+\kappa_{1}^{-1}}
=I_{t< T}
I_{s <t<s+\kappa_{1}^{-1}},
$$
 and that
\eqref{3.15.1} yields
$$
\int_{-\infty}^{\infty}I_{s\wedge T}(t)(\zeta^{s}_{t})^{p}\,ds
=\int_{-\infty}^{\infty}I_{s\wedge T<t<[(s\wedge T)+\kappa_{1}^{-1}]
\wedge T}
\zeta^{p}(t-s)\,ds
$$
$$
=I_{t<T}\int_{t-\kappa_{1}^{-1}}^{t}
\zeta^{p}(t-s)\,ds=I_{t<T}.
$$
 
 We also notice that, since $\kappa_{1}$
depends only on $d,p,\rho_{1},K$, we have
$$
\int_{-\infty}^{\infty}|\zeta'(s)|^{p}\,ds=N^{*}.
$$
 
Then we conclude   
$$
\lambda^{p/2} \| u\|^{p}_{\bL _{p}(T)} +
\| Du\|^{p}_{\bL _{p}(T)}
\leq N_{1} \sum_{i=1}^{d}\| 
  f^{i}\|^{p}_{\bL _{p}(T)}
+N_{1}\lambda^{-p/2 } \|  f^{0}\|^{p}_{\bL _{p}(T)}
$$
$$ 
+ N_{1}^{*} \lambda^{-p/2}\|    Du\|_{\bL_{p}(T)}^{p}
+ N_{1}^{*}\|   u\|_{\bL_{p}(T)}^{p}
+N_{1}^{*}\lambda^{-p/2}\sum_{i=1}^{d}\| 
 f^{i}\|^{p}_{\bL _{p}(T)}.
$$
Without losing generality we assume that $N_{1}\geq1$
and  we show how to choose $\lambda_{0}=\lambda_{0}(
d,\delta,p,\rho_{0},\rho_{1},K)\geq1$. 
Above we assumed that $\lambda\geq\lambda_{0}(d,\delta,p,\rho_{0})$,
where $\lambda_{0}(d,\delta,p,\rho_{0})$ is taken from
Lemma \ref{lemma 3.11.1}. Therefore, we take
$$
\lambda_{0}=\lambda_{0}(
d,\delta,p,\rho_{0},\rho_{1},K)\geq\lambda_{0}(d,\delta,p,\rho_{0})
$$
  such that
  $\lambda_{0}^{p/2}\geq 2N^{*}_{1} $. Then we obviously come
to \eqref{3.11.2} (with  $S=-\infty$). 
The theorem is proved.

\mysection{Proof of Theorems \protect\ref{theorem 3.16.1}
and \protect\ref{theorem 9.26.1}}
                                          \label{section 6.9.5}

We need two auxiliary results.

\begin{lemma}
                                         \label{lemma 9.2.1}
For any   $\tau,R\in(0,\infty)$,  we have  
\begin{equation}
                                                     \label{9.2.2}
\int_{-\tau}^{\tau}\int_{B_{R}}(|\gb_{s}(x)|^{p'}
+|b_{s}(x)|^{p'} +c_{s}^{p'}(x)) \,dxds<\infty.
\end{equation}
\end{lemma}

Proof. Obviously it suffices to prove \eqref{9.2.2}
with $B_{\rho_{1}}(x_{0})$ in place of $B_{R}$ for any $x_{0}$.
In that case, for instance,  (notice that $q\geq p'$,
see \eqref{8.11.3})
$$
\int_{B_{\rho_{1}}(x_{0})} |\gb_{s}(x) |^{p'}\,dx\leq N\big(
\int_{B_{\rho_{1}}(x_{0})} |\gb_{s}(x)-\bar{\gb}_{s}(x_{0})|^{q}\,dx
\big)^{p'/q}
+N|\bar{\gb}_{s}(x_{0})|^{p'}
$$
According to \eqref{6.11.2}
$$
\int_{B_{\rho_{1}}(x_{0})} |\gb_{s}(x) |^{p'}\,dx\leq N
+N|\bar{\gb}_{s}(x_{0})|^{p'}
$$
and in what concerns $\gb$ it only remains to use
Assumption \ref{assumption 2.7.2} (iii).
Similarly, $b_{s}$ and $c_{s}$ are treated. The lemma is proved.

The solution of our   equation will be obtained
as the weak limit of the solutions of equations
with cut-off coefficients. Therefore, the following result
is appropriate. By the way, observe that usual way
of proving the existence of solutions based on   a priori
estimates and the method of continuity
cannot work in our setting mainly because
of what is said after Theorem~\ref{theorem 8.8.3}.

\begin{lemma}
                                            \label{lemma 3.23.2}
Let $\phi\in C^{\infty}_{0}$, $\tau\in(0,\infty)$. 
Let $u^{m}$, $u\in \bW^{1}_{p}$, $m=1,2,...$, be such that
$u^{m}\to u$ weakly in $\bW^{1}_{p}$. 
For $m=1,2,...$   define $\chi_{m}(t)=(-m)\vee t\wedge m$,
$\gb^{i}_{mt}=\chi_{m}(\gb^{i}_{t})$,
$b^{i}_{mt}=\chi_{m}(b^{i}_{t})$, and $c_{mt}=\chi_{m}(c_{t})$.
Then  the functions
\begin{equation}
                                                 \label{4.19.6}
\int_{0}^{t} (b^{i}_{ms}D_{i}u^{m}_{s},\phi)\,ds,\quad
\int_{0}^{t} (\gb^{i}_{ms}u^{m}_{s},D_{i}\phi) \,ds,\quad
\int_{0}^{t} (c_{ms}u^{m}_{s}, \phi) \,ds
\end{equation}
converge weakly in the space $\cL_{p}([-\tau ,\tau ])$
as $m\to\infty$ to
\begin{equation}
                                                     \label{8.12.1}
\int_{0}^{t} (b^{i}_{s}D_{i}u_{s},\phi)\,ds,
\quad
\int_{0}^{t} (\gb^{i}_{s}u_{s},D_{i}\phi)\,ds,\quad
\int_{0}^{t} (c_{s} u_{s},\phi)\,ds,
\end{equation}
respectively.
\end{lemma}

Proof. By Corollary \ref{corollary 3.23.1} and by the
fact that (strongly) continuous operators are weakly
continuous we obtain that
$$
\int_{0}^{t} (b^{i}_{s}D_{i}u^{m}_{s},\phi) \,ds
\to
\int_{0}^{t} (b^{i}_{s}D_{i}u_{s},\phi) \,ds
$$
as $m\to\infty$ 
weakly in the space $\cL_{p}([-\tau ,\tau ])$.
 Therefore, in what concerns the first
function  in \eqref{4.19.6}, it suffices to show that
$$
\int_{0}^{t} (D_{i}u^{m}_{s},(b^{i}_{s}-b^{i}_{ms})\phi)
 \,ds
\to0
$$
weakly in $\cL_{p}([-\tau,\tau ])$. In other words, it suffices to show
that for any $\xi\in \cL_{p'}([-\tau ,\tau ])$
$$
\int_{-\tau}^{\tau}\xi_{t}
\big(\int_{0}^{t} (D_{i}u^{m}_{s},(b^{i}_{s}-b^{i}_{ms})\phi)
 \,ds
\big)\,dt\to0.
$$
This relation is rewritten as
\begin{equation}
                                                 \label{4.21.1}
 \int_{-\tau }^{\tau}
 (D_{i}u^{m}_{s},\eta_{s}(b^{i}_{s}-b^{i}_{ms})\phi)
 \,ds\to0,
\end{equation}
where  
$$
\eta_{s}:=\int_{s}^{\tau \text{\rm sgn}\, s }\xi_{t}\,dt
$$
is bounded on $[-\tau,\tau]$.
However, by  the 
dominated convergence theorem and Lemma \ref{lemma 9.2.1},
we have
$ 
\eta_{s}(b^{i}_{s}-b^{i}_{ms})\phi\to0 
$ 
as $m\to\infty$
strongly in $\bL_{p'}( -\tau,\tau  )$  and 
by assumption   $Du^{m}\to Du$
weakly in $\bL_{p}( -\tau ,\tau )$. This implies
\eqref{4.21.1}. Similarly,
one proves our assertion about
the remaining functions in \eqref{4.19.6}.  The lemma is proved.

{\bf Proof of Theorem \ref{theorem 3.16.1}}.
Owing to Theorem \ref{theorem 3.11.1} implying that 
the solution on $(-\infty,T]\cap\bR$ is unique,
without loss of generality we may assume that
$ T=\infty$. Define
$\gb_{mt}$, $b_{mt}$, and $c_{mt}$   as in  
 Lemma \ref{lemma 3.23.2} and consider equation
\eqref{2.6.4} with $\gb_{mt}$, $b_{mt}$, and $c_{mt}$  
in place of $\gb_{t}$, $b_{t}$, and $c_{t}$, respectively.
Obviously, $\gb_{mt}$, $b_{mt}$, and $c_{mt}$  satisfy Assumption
\ref{assumption 3.11.1} with the same $\gamma $ and $K$
as $\gb_{t}$, $b_{t}$, and $c_{t}$ do. By Theorem \ref{theorem 3.11.1}
and the method of continuity for  
$\lambda\geq\lambda_{0}(d,\delta,p,\rho_{0},\rho_{1},K)$ there exists
a unique solution $u^{m} $ of the modified
equation on $\bR$.

By Theorem \ref{theorem 3.11.1}  we also have
$$
\|u^{m}\|_{\bL_{p} }+\|Du^{m}\|_{\bL_{p} }\leq N,
$$
where $N$ is independent of $m$. Hence the sequence of functions
$u^{m} $ is bounded in the   space 
$\bW^{1}_{p}$ and consequently has a weak limit 
point $u\in \bW^{1}_{p}$. For simplicity of presentation
we assume that the whole sequence $u^{m}  $
converges weakly to $u$. 
Take a $\phi\in C^{\infty}_{0}$. Then
by Lemma \ref{lemma 3.23.2}   the functions \eqref{4.19.6}  
converge to \eqref{8.12.1} weakly in
$\cL_{p}([-\tau,\tau])$ as $m\to\infty$
for any $\tau$.  Obviously,
the same is true for $(u^{m}_{ t},\phi)\to(u_{t},\phi)$
and the remaining terms entering 
the equation for 
$ u^{m}_{t}$.
Hence, by passing to the weak limit in the equation
for $u^{m}_{ t}$ we see that for any
$\phi\in C^{\infty}_{0}$  equation
\eqref{3.16.7} holds for almost any $s,t\in\bR$. 

Now notice that, for each
$t\in\bR$, owing to Corollary \ref{corollary 6.27.1} the equation
$$
(\hat{u}_{t},\phi):=
\int_{0}^{1}(u_{s},\phi)\,ds+
\int_{0}^{1}\big( \int_{s}^{t}\big[(b^{i}_{r}D_{i}u_{r}
-(c_{r}+\lambda)u_{r}+f^{0}_{r},\phi)
$$
\begin{equation}
                                                    \label{8.11.1}
-(a^{ij}_{r}D_{j}u_{r}+\gb^{i}_{r}u_{r}+
f^{i}_{r},D_{i}\phi)
 \big]\,dr\big)\,ds
\end{equation}
defines a distribution. Furthermore, by the above
for any $\phi\in C^{\infty}_{0}$ we have
$(u_{t},\phi)=(\hat{u}_{t},\phi)$ (a.e.). A standard argument shows
that for almost all
$t\in\bR$,
$(u_{t},\phi)=(\hat{u}_{t},\phi)$ for any $\phi\in C^{\infty}_{0}$,
that is $u_{t}=\hat{u}_{t}$ (a.e.) and $\hat{u}_{t}\in\bW^{1}_{p}$.
In particular, we see that we can replace $u_{r}$ in \eqref{8.11.1}
with $\hat{u}_{r}$. Finally, for any $t_{1},t_{2}\in\bR$
$$
(\hat{u}_{t_{2}},\phi)-(\hat{u}_{t_{1}},\phi)=
\int_{0}^{1}\big(\int_{t_{1}}^{t_{2}}\big[(b^{i}_{r}D_{i}\hat{u}_{r}
-(c_{r}+\lambda)\hat{u}_{r}+f^{0}_{r},\phi)
$$
$$
-(a^{ij}_{r}D_{j}\hat{u}_{r}+\gb^{i}_{r}\hat{u}_{r}+
f^{i}_{r},D_{i}\phi)
 \big]\,dr\big)\,ds=\int_{t_{1}}^{t_{2}}\big[(b^{i}_{r}D_{i}\hat{u}_{r}
-(c_{r}+\lambda)\hat{u}_{r}+f^{0}_{r},\phi)
$$
 $$
-(a^{ij}_{r}D_{j}\hat{u}_{r}+\gb^{i}_{r}\hat{u}_{r}+
f^{i}_{r},D_{i}\phi)
 \big]\,dr
$$
and the theorem is proved.

{\bf Proof of Theorem \ref{theorem 9.26.1}}. (i) One reduces
the general case to the one that $u_{S}=0$ as in the proof
of Theorem \ref{theorem 8.8.1}. Also, obviously,
one can assume that $\lambda$ is as large as we like, say
satisfying \eqref{9.27.1},
since $S$ and $T$ are finite. By continuing $u_{t}(x)$ as zero
for $t\leq S$ we see that we may assume that $S=\infty$.
If we set $f^{j}_{t}=0$ for $t\geq T$ and use Theorem
\ref{theorem 3.16.1} about the existence of solutions
on $(-\infty,\infty)$
along with Theorem \ref{theorem 3.11.1},
which guarantees uniqueness of solutions on $(-\infty,T]$,
then we see that we only need to prove assertion (ii)
of the theorem.

(ii) In the above proof of Theorem \ref{theorem 3.16.1}
we have constructed the unique solutions of our
equations as the weak limits of the solutions of equations with
cut-off coefficients. Therefore, if we knew that the result
is true for equations with bounded coefficients, then we would obtain
it in our general case as well.

Thus it only remains to concentrate on equations
with bounded coefficients. Existence an uniqueness theorems
also show that it suffices to prove that, if $u$ is the solution
corresponding to $p=p_{2}$, then $u\in\bW^{1}_{p_{1}}$.

Take a   $\zeta\in C^{\infty}_{0}(\bR^{d+1})$ such that 
$\zeta(0)=1$,
set $\zeta^{n}_{t}(x)=\zeta(t/n,x/n)$,
and notice that $u^{ n}_{t}:=u_{t} \zeta^{n}_{t}$ satisfies
$$
\partial_{t}u^{ n}_{t}= L_{t}u^{ n}_{t} 
-\lambda u^{ n}_{t}+D_{i}f^{i}_{nt}+f^{ 0}_{nt}  ,
$$
where
$$
f^{i}_{nt}=f^{i}_{t}\zeta^{n}_{t}-u_{t}a^{ ij}_{t}
D_{j}\zeta^{n}_{t},\quad 
i\geq1,
$$
$$
f^{ 0}_{nt}=f^{0}_{t}\zeta^{n}_{t}-f^{i}_{t}D_{i}\zeta^{n}_{t}
-(a^{ij}_{t}D_{j}u_{t}+a^{i}_{t}u_{t})D_{i}\zeta^{n}_{t}
-b^{i}_{t}u_{t}D_{i}\zeta^{n}_{t}+u_{t}\partial_{t}\zeta^{n}_{t}.
$$
Since $u^{n}_{t}$ has compact support and $p_{1}\leq p_{2}$,
it holds that $u^{n}\in \bW^{1}_{p}$ for any $p\in [1,p_{2}]$
and by Theorem \ref{theorem 3.11.1}
  for $p\in [p_{1},p_{2}]$ we have
\begin{equation}
                                       \label{2.24.3}
\|u^{ n}\|_{\bW^{1}_{p} }
\leq N \sum_{i=0}^{d}\|f^{ i}_{n}\|_{\bL_{p} }.
\end{equation}
One knows that 
$$
\|f^{ i} \|_{\bL_{p} }\leq N(\|f^{ i} \|_{\bL_{p_{1} }}
+\|f^{ i} \|_{\bL_{p_{2}} }),
$$
so that by H\"older's inequality
$$
\|f^{ i}_{n}\|_{\bL_{p} }
\leq N+N\|uD\zeta^{n}\|_{\bL_{p} }
\leq N+ \|u\|_{\bL_{p_{2}}}\|D\zeta^{n}
\|_{\bL_{q}} ,
$$
with   constants  $N$ independent of $n$,
where
$$
q=\frac{pp_{2}}{p_{2}-p}.
$$
Similar estimates are available for other terms in
the right-hand side of \eqref{2.24.3}. Since
$$
\|\partial_{t}\zeta^{n}
\|_{\bL_{q} }+\|D\zeta^{n}
\|_{\bL_{q} }=Nn^{-1+(p_{2}-p)(d+1)/(p_{2}p)}\to0
$$
as $n\to\infty$ if
\begin{equation}
                                       \label{2.24.4}
\frac{1}{p}-\frac{1}{p_{2}}<\frac{1}{d+1},
\end{equation}
estimate \eqref{2.24.3} implies that $u 
\in\bW^{1}_{p} $.

Thus knowing that $u 
\in\bW^{1}_{p_{2}} $ allowed us to conclude
that $u 
\in\bW^{1}_{p} $ as long as $p\in[p_{1},p_{2}]$
and \eqref{2.24.4} holds. We can now replace $p_{2}$ 
with a smaller $p$ and keep going in the same way
each time increasing $1/p$ by the same amount until $p$
  reaches $p_{1}$. Then we get that 
$u 
\in\bW^{1}_{p_{1}} $. The theorem is proved

\end{document}